\documentclass[a4paper,11pt,twoside,]{article}
\usepackage[utf8]{inputenc}
\usepackage[T1]{fontenc}
\usepackage[sc]{mathpazo}
\usepackage[english]{babel}
\usepackage{geometry}
\usepackage{amsmath,amsfonts,amssymb,amsthm}
\usepackage[final]{hyperref}
  \usepackage{epsfig} 
  \usepackage{graphicx}  
\usepackage{epstopdf}
\usepackage{mathtools}
\usepackage[dvipsnames]{xcolor}
\usepackage{pdflscape}
\usepackage{etoolbox}
\usepackage{fancyhdr}
\usepackage{lastpage}
\usepackage{ifdraft}
\usepackage{hyphenat}
\ifdraft{\usepackage{showlabels}}{}
\usepackage{enumitem}
\usepackage{siunitx}
\usepackage[mathlines]{lineno}
\newcommand*\patchAmsMathEnvironmentForLineno[1]{%
  \expandafter\let\csname old#1\expandafter\endcsname\csname #1\endcsname
  \expandafter\let\csname oldend#1\expandafter\endcsname\csname end#1\endcsname
  \renewenvironment{#1}%
  {\linenomath\csname old#1\endcsname}%
  {\csname oldend#1\endcsname\endlinenomath}}% 
\newcommand*\patchBothAmsMathEnvironmentsForLineno[1]{%
  \patchAmsMathEnvironmentForLineno{#1}%
  \patchAmsMathEnvironmentForLineno{#1*}}%
\patchBothAmsMathEnvironmentsForLineno{equation}%
\patchBothAmsMathEnvironmentsForLineno{align}%
\patchBothAmsMathEnvironmentsForLineno{flalign}%
\patchBothAmsMathEnvironmentsForLineno{alignat}%
\patchBothAmsMathEnvironmentsForLineno{gather}%
\patchBothAmsMathEnvironmentsForLineno{multline}%
\newcommand{\dd}{\mathop{}\!\mathrm{d}}
\theoremstyle{plain}

\numberwithin{equation}{section}
\numberwithin{table}{section}
\numberwithin{figure}{section}

\title{Efficient numerical computation of the basic reproduction number for structured populations}
\author{
Dimitri Breda$^{1}$, Francesco Florian$^{2}$, Jordi Ripoll$^{3}$ and Rossana Vermiglio$^{4}$\\[.5em]
\small $^{1}$CDLab -- Computational Dynamics Laboratory\\
\small Department of Mathematics, Computer Science and Physics -- University of Udine\\
\small via delle scienze 206, 33100 Udine, Italy\\
\small dimitri.breda@uniud.it\\[.3em]
\small $^{2}$ Institute of Mathematics--University of Z\"{u}rich\\
\small Winterthurerstrasse 190, CH-8057 Z\"{u}rich, Switzerland\\
\small and CDLab -- Computational Dynamics Laboratory\\
\small francesco.florian@uzh.ch\\[0.3em]
\small $^{3}$Department of Computer Science, Applied Mathematics and Statistics -- University of Girona\\
\small Campus Montilivi, 17003 Girona, Spain\\
\small jripoll@imae.udg.edu\\[.3em]
\small $^{4}$CDLab -- Computational Dynamics Laboratory\\
\small Department of Mathematics, Computer Science and Physics -- University of Udine\\
\small via delle scienze 206, 33100 Udine, Italy\\
\small rossana.vermiglio@uniud.it\\[-0.5cm]
}
\date{}
\begin{document}

\maketitle

\begin{abstract} As widely known, the basic reproduction number plays a key role in weighing birth/infection and death/recovery processes in several models of population dynamics. In this general setting, its characterization as the spectral radius of next generation operators is rather elegant, but simultaneously poses serious obstacles to its practical determination. In this work we address the problem numerically by reducing the relevant operators to matrices through a pseudospectral collocation, eventually computing the sought quantity by solving finite-dimensional eigenvalue problems. The approach is illustrated for two classes of models, respectively from ecology and epidemiology. Several numerical tests demonstrate experimentally important features of the method, like fast convergence and influence of the smoothness of the models' coefficients. Examples of robust analysis of instances of specific models are also presented to show potentialities and ease of application. 


\end{abstract}
%\communicated{...}
%%\dedication{...}
%\received{...}
%\accepted{...}
%\journalname{...}
%\journalyear{...}
%\journalvolume{..}
%\journalissue{..}
%\startpage{1}
%\aop
%\DOI{...}

%-----------------------------------------------------------------------------
%-----------------------------------------------------------------------------
\section{Introduction}
\label{s_introduction}
Aim of this paper is to give a concise, efficient and effective method for the numerical computation of the {\it basic reproduction number}, a fundamental tool in the analysis of ecological and epidemiological models of population dynamics, as well as in tackling evolutionary aspects in those populations, see e.g. \cite[Section 2.1]{bcr17b} and the references therein. We refer to \cite{hes02} as a starting reference on this quantity, commonly denoted $R_{0}$.

\bigskip
The need for a numerical approach is soon justified since we are interested in population models posed on infinite-dimensional spaces, where $R_{0}$ is usually characterized as the spectral radius of suitable operators, mainly known as {\it next generation operators}, see, e.g., \cite{dhm90}. As in general there is no explicit expression for $R_{0}$ readily practicable for computation (and in some cases nor the relevant operators are explicitly known), the main idea is to reduce to finite dimension through some discretization techniques. Then approximations to $R_{0}$ are obtained by computing the (dominant) eigenvalues of the resulting matrices.

As far as this general strategy is concerned, to the best of our knowledge \cite{kun17} is the first scheme appearing in the literature. Therein an elementary discretization based on the celebrated Euler's method is proposed to approximate the basic reproduction number of a class of age-structured epidemics. Very recently a slight improvement has been presented in \cite{guo19}, obtained by employing $\theta$-methods, always in the context of age-structured epidemics. Made exception for \cite{fv19}, which is anyway related to the approach proposed here (see below), \cite{guo19,kun17} are the only two works concerning the proper numerical computation of $R_{0}$ in continuously structured populations.

Here, on the one hand, we intend to largely improve the seminal idea of \cite{kun17} by applying state-of-the-art discretizations based on pseudospectral polynomial collocation (see, e.g., \cite{ab19,and19}). The general outcome is a more reliable tool, with faster convergence ideally of infinite order (known as {\it spectral accuracy}, see, e.g., \cite{tref00}). This represents a great advantage compared to the finite-order convergence of \cite{guo19,kun17}, as it translates into much more accurate approximations obtained with much smaller matrices. As a consequence the required computational load reduces, in terms of both time and memory usage. The latter is a favorable feature when stability and bifurcation analyses are the final target in presence of varying or uncertain model parameters, and of course this is frequently the case in realistic contexts.

On the other hand, we aim at increasing the flexibility of the numerical approach, widening the applicability also to other models, e.g., from ecology. In particular, we illustrate two paradigmatic classes of models, namely spatially-structured cell populations (whose prototype representative is referred to as {\it model A} in the sequel) and infectious diseases in age-structured populations (whose prototype representative is referred to as {\it model B} in the sequel). This choice potentially serves as a basis for extending the proposed method to more involved systems of population dynamics (e.g., with several structuring variables). If a structured epidemic model describing the novel pandemia of coronavirus was available, such as {\it Model B} or any of its extensions, 
this method could be applied to compute $R_0$ for different values of the parameters.

\bigskip
The contents of this research are organized as follows. In Section \ref{s_models} we present the two classes of models mentioned above, introducing all the necessary ingredients and the operators leading to the definition of $R_{0}$. The numerical treatment is proposed in Section \ref{s_discretization}, where explicit expressions for the discretizing matrices are constructively given, also for the sake of implementation for those interested (Matlab codes are freely available at \url{http://cdlab.uniud.it/software}). In Section \ref{s_results} we first illustrate the numerical properties of the proposed techniques, analyzing experimentally error and convergence, and then we finally perform some robust analysis on specific instances of the models of interest.

\bigskip
Let us remark that here we deliberately quit to tackle a theoretical investigation of the convergence features of the proposed methods. This would indeed require suitable tools, falling out of the scope of the current experimental study, as well as enough room for a proper description. We thus reserve to give formal proofs as well as possible improvements in a related work, currently in preparation \cite{brv20}. As far as other developments are concerned, we also plan to extend this approach in order to treat periodic problems, that is, to population models where the environment is time-periodically varying. In this respect we suggest the reading of the very recent work \cite{ina19}.

\bigskip
Finally, we wish to mention that the present research originated from \cite{flo18}, where an embryonic study of collocation techniques for the numerical computation of $R_{0}$ is proposed (also the work \cite{fv19} mentioned above is inspired from \cite{flo18}). Moreover, let us also note that pseudospectral techniques have become a reference tool in the numerical treatment of infinite-dimensional dynamical systems in view of either stability and bifurcation analyses in the context of population dynamics, see, e.g., \cite{bdgsv16} or \cite{bmv15} and the references therein.
%-----------------------------------------------------------------------------
%-----------------------------------------------------------------------------
\section{Models and theoretical background}
\label{s_models}

Continuously structured  models of population dynamics can be described by nonlinear abstract differential equations for the density of individuals with respect to some structuring variables, e.g., age, size or space.

With the aim at setting the theoretical background, let $X$ be a Banach lattice and let $u\in X$ represent the density of individuals. We introduce next a continuous-time evolution equation for the population density. For the theory of population dynamics we refer, among many others, to the monographs \cite{dhb13,ian95,im17,ip14,ina17,thi03}.

\bigskip
On the one hand, if we focus on {\it ecological} models (the class of models A in the sequel), then in general a nonlinear problem can be decomposed into birth and non-birth terms as 
\begin{equation}\label{nonlinear}
u'(t)=B(u(t))u(t)-M(u(t))u(t),\quad u(t)\in X,\;t\geq0,
\end{equation}
where, for any $u\in X$, $B(u)$ is the birth linear operator and $M(u)$ is the linear operator accounting for other processes non related to birth (e.g., mortality or transition). Of course, this decomposition depends on what is exactly considered as a birth event, see an example in Section \ref{s_modelA} and also \cite{bcr17a,cd16}.

If we are concerned with the extinction of the population in \eqref{nonlinear}, then we can focus on the steady-state $u^{\ast}\equiv0$. Its stability is analyzed by means of the (formally) linearized model
\begin{equation}\label{1bis}
u'(t)=B(0)u(t)-M(0)u(t),\quad u(t)\in X,\;t\geq0.
\end{equation}

\bigskip
On the other hand, if we focus on {\it epidemiological} models (the class of models B in the sequel), then we firstly need to distinguish between infective individuals $u\in X$ and (multi-type set of) non-infective individuals $v\in X^{m}$, with $m$ some positive integer. Here, the original nonlinear problem consists of $m+1$ equations for the population densities of the two types of individuals. The one relevant to the infective individuals can be decomposed according to infection and non-infection terms as
\begin{equation}\label{nonlinear2a}
u'(t)=B(v(t),u(t))u(t)-M(v(t),u(t))u(t),\quad v(t)\in X^{m},\;u(t)\in X,\;t\geq0.
\end{equation}
The remaining ones have the general form
\begin{equation}\label{nonlinear2b}
v'(t)=F(v(t),u(t)),\quad v(t)\in X^{m},\;u(t)\in X,\;t\geq0,
\end{equation}
with $F$ depending on the type of epidemics. For any $v\in X^{m}$ and $u\in X$, $B(v,u)$ is the linear infection operator and $M(v,u)$ is the linear operator describing processes other than infection (e.g., permanent or temporal recovery or transition processes). Again, the decomposition may differ depending on what is considered as an infection event. 

Here, typically, the main concern is the early stage of the infection, i.e., the initial epidemic spread in \eqref{nonlinear2a}-\eqref{nonlinear2b}. We can thus focus on the disease-free steady state $(v^{\ast},u^{\ast})\in X_{+}^{m}\times\{0\}$ (where $X_{+}$ stands for the positive cone in $X$). Since the (formally) linearized system has triangular form, under common assumptions, the stability analysis of the disease-free equilibrium is reduced to the linear model\begin{equation}\label{3bis}
u'(t)=B(v^{\ast},0)u(t)-M(v^{\ast},0)u(t),\quad u(t)\in X,\;t\geq0.
\end{equation}

\bigskip
Summarizing and simplifying the notation to make it uniform in \eqref{1bis} and \eqref{3bis}, for either model A or B we are faced with the analysis of an abstract linear equation of the form
\begin{equation}\label{Model}
u'(t)=Bu(t)-Mu(t),\quad u(t)\in X,\;t\geq0.
\end{equation}
We recall that above $X$ is a Banach lattice where the population density lives, $B:X\to X$ is a linear operator accounting for the {\it birth} process (meant as proper birth, or infection) and $M:\mathcal{D}(M)\subseteq X\to X$ is a linear operator accounting for the remaining processes, which we call {\it mortality} for brevity (thus meant as proper mortality, or recovery or, in general, any stage transition as in \cite[Part 2]{thi03}). Typically, the domain $\mathcal{D}(M)$ of $M$ is a subset of a subspace $Y\subseteq X$ where some degree of smoothness is present, characterized by additional constraints in general described through a linear map $\mathcal{C}:X\to\mathbb{R}^{p}$ for some positive integer $p$:
\begin{equation}\label{DM}
\mathcal{D}(M)=\{\phi\in Y\ :\ \mathcal{C}\phi=0\}.
\end{equation}
Other common (biologically meaningful) assumptions are required. In particular, as in \cite{bcr17b}, $B$ is positive and bounded, while $-M$ is the infinitesimal generator of a strongly-continuous semigroup $\{T(t)\}_{t\geq0}$ of positive linear operators. Its spectral bound $s(-M)$ is strictly negative, which accounts for extinction in absence of birth. Consequently, $0$ belongs to the resolvent of $M$ and thus the latter is invertible with $M^{-1}=\int_{0}^\infty T(t)\dd t$. For these and other aspects on (positive) linear operators see, e.g., \cite{sha74,nee96}.

\bigskip
It is worth to mention that although \eqref{Model} represents a rather large family of population models (including in particular those with discrete structure when $\dim X<+\infty$), see, e.g., \cite[Section 7.2]{dhb13}), not all the structured models can be cast into this form, see the discussion in \cite{bcr17b}. However, in this standard case we can define the next generation operator as the operator $BM^{-1}:X\to X$ acting on the same (whole) space of the population density and, eventually, characterize the basic reproduction number as its spectral radius:
\begin{equation}\label{R0}
R_{0}:=\rho(BM^{-1}).
\end{equation}
Let us remark that the next generation operator is a generalization to the infinite-dimensional setting of the well-known next generation matrix for finite-dimensional models (i.e., $\dim X <+\infty$). Moreover, although the theoretical framework for the basic reproduction number has been well-established by many authors, see \cite{dhb13,dhm90,hes02,ina17} and the references therein to name a few, there is a lot of room for its efficient computation as already anticipated in the introduction. 

It follows from definition \eqref{R0} that $R_{0}$ is a non-negative spectral value being $BM^{-1}$ positive and bounded, see, e.g., \cite{sha74}. It is actually the largest $\lambda\geq0$ for which $B\phi-\lambda M\phi=\xi$ cannot be solved uniquely for $\phi\in\mathcal{D}(M)$ once $\xi\in X$ is given, see, e.g., \cite{bcr17b}. Moreover, if it happens that this $\lambda$ is a generalized eigenvalue, then
\begin{equation}\label{geig}
B\phi=\lambda M\phi
\end{equation}
for some nontrivial generalized eigenfunction $\phi\in\mathcal{D}(M)$. %or, equivalently,
%\begin{equation}\label{eig}
%BM^{-1}\psi=\lambda\psi
%\end{equation}
%for some nontrivial eigenfunction $\psi$ ($=M\phi$).
In particular, if $BM^{-1}$ is also compact, then the Krein-Rutman theorem \cite{kr48} ensures that $R_{0}$ is a positive eigenvalue (for slightly weaker, yet involved assumptions see \cite[Section 10.2]{ina17}). At this point let us anticipate that the numerical method proposed in Section \ref{s_discretization} relies on this assumption, so that \eqref{geig} can be reduced to a standard (read finite-dimensional) generalized eigenvalue problem for matrices. Compactness is proved in \cite{brv20} under reasonable assumptions on the specific ingredients of models A and B as described in the forthcoming Sections \ref{s_modelA} and \ref{s_modelB}. A discussion on further aspects related to $R_{0}$ being or not a generalized eigenvalue and their relevance on the numerical outcome is left to Section \ref{s_results}.

Following definition \eqref{R0}, with any further assumption, we obtain the useful upper bound
\begin{equation}\label{upperR0}
R_{0}\leq\|BM^{-1}\|_{X\leftarrow X}:=\sup_{\psi\in X\setminus\{0\}}\frac{\|BM^{-1}\psi\|_{X}}{\|\psi\|_{X}}=\sup_{\phi\in\mathcal{D}(M)\setminus\{0\}}\frac{\|B\phi\|_{X}}{\|M\phi\|_{X}}.
\end{equation}
If, in addition, $R_{0}=\lambda\geq0$ is a generalized eigenvalue with generalized eigenfunction $\phi\in\mathcal{D}(M)\setminus\{0\}$, then
\begin{equation}\label{upperR0eig}
R_{0}=\frac{\|B\phi\|_{X}}{\|M\phi\|_{X}}.
\end{equation}

\bigskip
Finally, let us also recall that $R_{0}$ is an alternative to the Malthusian parameter, i.e., the exponential population growth rate. Concerning computations, in general it can be more difficult to deal with the latter%Malthusian parameter of \eqref{Model}
, i.e., the spectral bound $r:=s(B-M)$ of the complete generator, rather than with the basic reproduction number $R_{0}$. For instance, we cannot ensure that $r$ is a spectral value (not even an eigenvalue) and there can be no upper bound for $r$ as it is always the case for $R_{0}$ (corresponding to a sufficient condition for population extinction or infection eradication). Moreover, the rank of $B$ is typically lower than that of $B-M$, which is an advantage for $R_{0}$. Finally, extra evolutionary aspects added to the models are more often related to $R_{0}$ than to $r$. Summarizing, besides the well-known sign relation between these quantities, see, e.g., \cite{thi09}, the basic reproduction number offers diverse advantages over the Malthusian parameter, from both the theoretical and the numerical points of view.
%-----------------------------------------------------------------------------
\subsection{Model A: spatially-structured cell populations}
\label{s_modelA}
As a first representative of prototype models to illustrate our numerical approach, let us consider a population of bacteria proliferating and moving along the intestine of an animal host, see, e.g., \cite{bcr17a,bcr17b}. Given that the intestine is portrayed as the line segment $[0,l]$ of finite length $l>0$, let us set $X:=L^{1}(0,l)$ and let $u(\cdot,t)\in X$ be the spatial density of bacteria along the intestine at time $t\geq0$. The nonlinear problem is described as
\begin{equation}\label{cell}
\left\{\setlength\arraycolsep{0.1em}\begin{array}{l}
\partial_{t} u(x,t)+\partial_{x}\big[c(x,S(t))u(x,t)-D(x,S(t))\partial_{x} u(x,t)\big]\\[1mm]
\qquad+\mu(x,S(t))u(x,t)=\beta(x,S(t))u(x,t)\\[1mm]
c(0,S(t))u(0,t)-D(0,S(t))\partial_{x} u(0,t)=0\\[1mm]
c(l,S(t))u(l,t)-D(l,S(t))\partial_{x} u(l,t)=0,
\end{array}\right.
\end{equation}
where $S(t):=\int_{0}^{l}\sigma(x)u(x,t)\dd x$, $t\geq0,$
%\begin{equation*}
%S(t):=\int_{0}^{l}\sigma(x)u(x,t)\dd x,\quad t\geq0,
%\end{equation*}
is the population size weighted through a non-negative weight $\sigma$. Above, $c$ is the velocity of the flow, $D$ is the diffusion coefficient, $\beta$ is the fertility rate and $\mu$ is the mortality rate. All are non-negative functions of the position $x\in[0,l]$ and of the size $S$. Note that transport, diffusion and vital processes are density-dependent accounting for limited resources (crowding effects) and space-specific accounting for the heterogeneity of the intestine, see \cite{bcr17a,bcr17b,ip14}.

\bigskip
For a population of bacteria, one can consider either symmetric division, i.e., when a mother cell divides then two daughter cells are born and the former disappears, or asymmetric division, i.e., one mother and one daughter. Both types can be considered simultaneously in \eqref{cell} by setting $\beta-\mu=\big[\nu\beta+(1-\nu)2\beta\big]-\big[(1-\nu)\beta+\mu\big]$
%\begin{equation*}
%\beta-\mu=\big[\nu\beta+(1-\nu)2\beta\big]-\big[(1-\nu)\beta+\mu\big]
%\end{equation*}
for some $\nu\in[0,1]$ representing the probability of asymmetric division. For instance, for space-independent vital rates one would readily get that $R_{0}=(2-\nu)\beta/[(1-\nu)\beta+\mu]$ since integrating over the whole space in \eqref{cell}, the system reduces to an ODE for the population size, see \cite{bcr17b}. The latter is a clear example of the fact that the basic reproduction number depends on what is exactly considered as a birth event. In any case, here we assume symmetric division ($\nu=0$) for simplicity, and without loss of generality with respect to the numerical discretization.

\bigskip
In what follows we are concerned with the extinction of the bacterial population. Linearizing \eqref{cell} around the extinction equilibrium $u^{\ast}\equiv0$ gives
\begin{equation}\label{linear-cell}
\left\{\setlength\arraycolsep{0.1em}\begin{array}{l}
\partial_{t} u(x,t)+\partial_{x}\big[c(x)u(x,t)-D(x)\partial_{x} u(x,t)\big]\\[1mm]
\qquad+[\beta(x)+\mu(x)]u(x,t)=2\beta(x)u(x,t)\\[1mm]
 c(0)u(0,t)-D(0)\partial_{x}u(0,t)=0\\[1mm]
 c(l)u(l,t)-D(l)\partial_{x}u(l,t)=0,
\end{array}\right.
\end{equation}
where we implicitly set $c(x):=c(x,0)$, $D(x):=D(x,0)$, $\beta(x):=\beta(x,0)$ and $\mu(x):=\mu(x,0)$ to lighten the notation though with a little abuse. Moreover, we further assume $c(x)>0$ for all $x\in(0,l)$, as well as $\beta$ and $\mu$ bounded with their sum not identically vanishing. Concerning diffusion, we consider either $D(x)\geq\tilde{D}$, $x\in[0,l]$, for some positive $\tilde{D}$ (everywhere positive diffusion) or $D\equiv0$ (complete lack of diffusion).

\bigskip
Eventually, with respect to the reference linear equation \eqref{Model}, it is natural to define from \eqref{linear-cell} the birth operator $B:X\rightarrow X$ as
\begin{equation}\label{BA}
(B\phi)(x):=2\beta(x)\phi(x),\quad x\in[0,l],
\end{equation}
and the mortality operator $M:\mathcal{D}(M)\subseteq X\rightarrow X$ as
\begin{equation}\label{MA}
(M\phi)(x):=\left[c(x)\phi(x)-D(x)\phi'(x)\right]'+[\beta(x)+\mu(x)]\phi(x),\quad x\in[0,l],
\end{equation}
with domain
\begin{equation}\label{DMA}
\setlength\arraycolsep{0.1em}\begin{array}{rcl}
\mathcal{D}(M):=\big\{\phi\in X\ &:\ & \phi',(c\phi-D\phi')'\in X\text{ and }\\[1mm]
&&c(x)\phi(x)-D(x)\phi'(x)=0\text{ for } x=0,l\big\}.
\end{array}
\end{equation}
Then, following \eqref{upperR0}, the upper bound
\begin{equation}\label{R0<2}
R_{0}\leq\sup_{\phi\in\mathcal{D}(M)\setminus\{0\}}\frac{2\int_{0}^{l}\beta(x)|\phi(x)|\dd x}{\int_{0}^{l}[\beta(x)+\mu(x)]|\phi(x)|\dd x}\leq2
\end{equation}
holds, as expected due to symmetric division.%-----------------------------------------------------------------------------
\subsection{Model B: age-structured epidemics}
\label{s_modelB}
As a second representative of prototype models to illustrate our numerical approach, let us consider the spread of an epidemics in an age-structured population. We split up the individuals of the population into three classes according to the stage of the infectious disease. If $l>0$ denotes the maximum (chronological) age, we set $X:=L^{1}(0,l)$ and let $s(\cdot,t)$, $i(\cdot,t)$ and $r(\cdot,t)$ in $X$ denote the density with respect to age of, respectively, susceptible, infected and removed individuals at time $t\geq0$, see, e.g., \cite{ian95,ina17}. The nonlinear problem is described as
\begin{equation}\label{age-epidemic}
\left\{\setlength\arraycolsep{0.1em}\begin{array}{l}
\partial_{t} s(x,t)+\partial_{x}s(x,t)+\mu(x)s(x,t)=-f(x,t)s(x,t)+\delta(x)i(x,t)+\sigma(x)r(x,t)\\[1mm]
\partial_{t} i(x,t)+\partial_{x}i(x,t)+\mu(x)i(x,t)=f(x,t)s(x,t)-[\gamma(x)+\delta(x)]i(x,t)\smallskip\\
\partial_{t} r(x,t)+\partial_{x}r(x,t)+\mu(x)r(x,t)=\gamma(x)i(x,t)-\sigma(x)r(x,t)\\[1mm] s(0,t)=\int_{0}^{l}\beta_{0}(x)[s(x,t)+(1-\theta)i(x,t)+(1-\hat{\theta})r(x,t)]\dd x\\[1mm] i(0,t)=\theta\int_{0}^{l}\beta_{0}(x)i(x,t)\dd x\\[1mm]
r(0,t)=\hat{\theta}\int_{0}^{l}\beta_{0}(x)r(x,t))\dd x,
\end{array}
\right.
\end{equation}
where $f$ is the force of infection defined as $f(x,t):=\int_{0}^{l}K_{0}(x,y,N(t))i(y,t)\dd y$, $x\in[0,l]$, $t\geq0$, with $N(t)$ the total population at time $t$ and infection kernel $K_{0}$ typically decomposed as $K_{0}(x,y,N)=\chi(y)\cdot c(x,y,N)/N$, $x,y\in[0,l]$, $N>0$. Above, $\chi$ is the probability of infection transmission through an infectious contact and $c$ is the density-dependent contact rate of a susceptible individual of age $x$ with an infected individual of age $y$. Moreover, $\beta_{0}$, $\mu$, $\gamma$, $\delta$ and $\sigma$ are, respectively, the fertility, (natural) mortality, removal, recovery and loss of immunity rates, all non-negative and with $\beta_{0}$ not identically vanishing. Finally, $\theta$ and $\hat{\theta}$ are the probabilities of vertical transmission of infectiveness and immunity, respectively. Note that recovery, removal, immunity loss and vital processes are all age-specific during the lifespan $[0,l]$.

\bigskip
System \eqref{age-epidemic} includes different types of epidemic models. For instance, the SIS model is given by $\gamma\equiv\sigma\equiv0$; the SIR model is given by $\delta\equiv\sigma\equiv0$; the SIRS model is given by $\delta\equiv0$. The exact interpretation of the removed class depends on the type of model considered. Moreover, we could also consider a vaccination rate and the resulting system would have an extra transition, from susceptible to recovered individuals in this case, or the consideration of exposed individuals, or even we could consider a multistrain epidemic model. See \cite{km33} and \cite[Chapter VI]{ian95}, \cite[Chapters 6--8]{ina17}, \cite[Chapter 22]{thi03}.

\bigskip
From now on, we assume that the population has reached a demographic steady state. Indeed, we assume zero demographic growth, i.e.,
\begin{equation}\label{intbeta0Pi=1}
\int_{0}^{l}\beta_{0}(x)\Pi_{0}(x)\dd x=1,
\end{equation}
where $\Pi_{0}(x):=\exp\left(-\int_{0}^{x}\mu(y)\dd y\right)$, $x\in[0,l]$, together with $\Pi_{0}(l)=0$, defines the survival probability, as well as initial conditions such that $s_{0}(x)+i_{0}(x)+r_{0}(x)=N_{0}\Pi_{0}(x)/\int_{0}^{l}\Pi_{0}(y)\dd y$, $x\in[0,l]$, for some positive $N_{0}$. Therefore, the total population $N(t)=\int_{0}^{l}[s(x,t)+i(x,t)+r(x,t)]\dd x=N_{0}$, $t\geq0$, remains constant over time.

\bigskip
In what follows we are concerned with the early stage of the epidemics of an infectious disease of any type (SIS, SIR or SIRS). Linearizing \eqref{age-epidemic} around the disease-free equilibrium $s^{\ast}(x)=N_{0}\Pi_{0}(x)/\int_{0}^{l}\Pi_{0}(y)\dd y$, $i^{\ast}\equiv0$, $r^{\ast}\equiv0$ gives
\begin{equation*}
\left\{\setlength\arraycolsep{0.1em}\begin{array}{l}
\partial_{t} i(x,t)+\partial_{x}i(x,t)=\int_{0}^{l}\chi(y)\cdot\frac{c(x,y,N_{0})}{N_{0}}\cdot i(y,t)\dd y\cdot\frac{N_{0}\,\Pi_{0}(x)}{\int_{0}^{l}\Pi_{0}(y)\dd y}\\[1mm]
\qquad-[\mu(x)+\gamma(x)+\delta(x)]i(x,t)\\[1mm]
i(0,t)=\theta\int_{0}^{l}\beta_{0}(x)i(x,t)\dd x
\end{array}\right.
\end{equation*}
for the infected individuals. We let $\int_{0}^{l}\mu(x)\dd x=+\infty$ to avoid immortal individuals. In order to properly deal with this condition, we make the change of variables $i(x,t)=u(x,t)\Pi_{0}(x)$, thus removing the mortality term and reducing the problem to
\begin{equation}\label{linear-epi}
\left\{\setlength\arraycolsep{0.1em}\begin{array}{l}
\partial_{t} u(x,t)+\partial_{x}u(x,t)=\int_{0}^{l}K(x,y)u(y,t)\dd y-[\gamma(x)+\delta(x)]u(x,t)\\[1mm]
u(0,t)=\theta\int_{0}^{l}\beta(x)u(x,t)\dd x 
\end{array}\right.
\end{equation}
with effective infection kernel 
\begin{equation}\label{effective-kernel}
K(x,y):=\chi(y)c(x,y,N_{0})\cdot\frac{\Pi_{0}(y)}{\int_{0}^{l}\Pi_{0}(z)\dd z}\geq0,\quad x,y\in[0,l],
\end{equation}
and effective fertility rate $\beta(x):=\beta_{0}(x)\Pi_{0}(x)\geq0$, $x\in[0,l]$. We recall that $\int_{0}^{l}\beta(x)\dd x=1$ due to \eqref{intbeta0Pi=1}. Note that here $u(l,t)$ can be strictly positive allowing for the possibility of chronic infected individuals.

\bigskip
Eventually, with respect to the reference linear equation \eqref{Model}, from the normalized model \eqref{linear-epi} it is natural to define the birth\footnote{We recall that it is actually the infection operator.} operator $B:X\rightarrow X$ as
\begin{equation}\label{BB}
(B\phi)(x):=\int_{0}^{l}K(x,y)\phi(y)\dd y,\quad x\in[0,l],
\end{equation}
and the mortality\footnote{We recall that it is actually the recovery and transition operator.} operator $M:\mathcal{D}(M)\subseteq X\rightarrow X$ as
\begin{equation}\label{MB}
(M\phi)(x):=\phi'(x)+[\gamma(x)+\delta(x)]\phi(x),\quad x\in[0,l],
\end{equation}
with domain
\begin{equation}\label{DMB}
\mathcal{D}(M):=\left\{\phi\in X\ :\ \phi'\in X\text{ and }\phi(0)=\theta\int_{0}^{l}\beta(x)\phi(x)\dd x\right\}.
\end{equation}

Moreover, we can explicitly compute the (normalized) next generation operator $BM^{-1}$ (and even obtain explicit expressions for some special cases, see Section \ref{s_results}). Indeed, by defining $\Pi_{1}(x):=\exp\left( -\int_{0}^{x}[\gamma(y)+\delta(y)]\dd y\right)\leq 1$, the next generation operator becomes
\begin{equation}\label{NG-Age-Epidemic}
(BM^{-1}\psi)(x)=\int_{0}^{l}K(x,y)\Pi_{1}(y)\left(C+\int_{0}^{y}\frac{\psi(z)}{\Pi_{1}(z)}\dd z\right)\dd y,\quad x\in[0,l],
\end{equation}
with 
\begin{equation*}
C:=\theta\int_{0}^{l}\beta(x)\left(\int_{0}^{x}\frac{\Pi_{1}(x)}{\Pi_{1}(z)}\psi(z)\dd z\right)\dd x\cdot\left(1-\theta\int_{0}^{l}\beta(x)\Pi_{1}(x)\dd x\right)^{-1}\geq0.
\end{equation*}
Above, we implicitly assume that $\theta\Pi_{1}\not\equiv 1$ in order for the last factor to be well-defined. Finally, following \eqref{upperR0}, we obtain the upper bound
\begin{equation}\label{upperR0B}
R_0\leq \int_{0}^{l}\int_{0}^{l}K(x,y) \left( {\textstyle\frac{\theta \Pi_{1}(y)}{1-\theta\int_{0}^{l}\beta(x)\Pi_{1}(x)\dd x} + 1} \right)   \dd y\dd x,
\end{equation}
%\begin{equation}\label{upperR0B}
%R_{0}\leq\frac{\int_{0}^{l}\int_{0}^{l}K(x,y)\dd y\dd x}{1-\theta\int_{0}^{l}\beta(x)\Pi_{1}(x)\dd x},
%\end{equation}
representing the total number of infections taking into account vertical transmission. Eventually, the rhs of \eqref{upperR0B} being less than 1, can be used as a sufficient condition for the disease eradication.
%-----------------------------------------------------------------------------
%-----------------------------------------------------------------------------
\section{Discretization}
\label{s_discretization}
Let $N$ be a positive integer. Under the assumption that $R_{0}$ is a generalized eigenvalue, we propose a general discretization approach consisting in the collocation of \eqref{geig} on a mesh of $N+1$ distinct nodes $0=:x_{N,0}<x_{N,1}<\cdots<x_{N,N}:=l$ distributed on $[0,l]$, $l>0$. Let us soon observe that collocation is meaningless in the whole space $X=L^{1}(0,l)$, but the method is restricted to generalized eigenfunctions $\phi$, which are, in general, smooth enough to guarantee pointwise definition, recall indeed \eqref{DM}.

\bigskip
In the sequel, let $X_{N}:=\mathbb{R}^{N+1}$ be the finite-dimensional counterpart of $X$ and let $\Phi:=(\Phi_{0},\Phi_{1},\ldots,\Phi_{N})^{T}\in X_{N}$ with $\Phi_{i}$ representing the numerical approximation of $\phi(x_{N,i})$, $i=0,1,\ldots,N$.

The choice of including $x_{N,0}=0$ and $x_{N,N}=l$ in the collocation mesh is justified by the fact that, in general, $\phi$ belongs to some restricted domain \eqref{DM} where the mortality operator $M$ is suitably defined. For the models of interest in this work, in fact, this domain is characterized by boundary conditions at one or both the extrema $0,l$ of the domain of the function space $X$, see indeed \eqref{DMA} and \eqref{DMB}. Therefore, the choice above allows for a simplified treatment as one can directly relate the boundary conditions to the collocation equations. 

As the concerned mortality operators usually involve differentiation and/or integration, see \eqref{MA} and \eqref{MB}, let us also introduce the {\it differentiation matrix} and the {\it quadrature weights} associated to the collocation nodes. The first, denoted $H_{N}$, has entries
\begin{equation*}
h_{N;i,j}:=\ell_{N,j}'(x_{N,i}),\quad i,j=0,1,\ldots,N,\footnote{Here $a_{N;i,j}$ denotes the $(i,j)$-th entry of a matrix $A\in\mathbb{R}^{(N+1)\times(N+1)}$, $i,j=0,1,\ldots,N$. This notation is also extended later on according to Matlab's {\it colon} notation: for instance $a_{N;i,0:N}$ and $a_{N;0:N,j}$ indicate, respectively, the $i$-th row and the $j$-th column of $A$. }
\end{equation*}
where $\{\ell_{N,0},\ell_{N,1},\ldots,\ell_{N,N}\}$ is the Lagrange basis relevant to the chosen nodes: if $f$ is a smooth function on $[0,l]$ and $v:=(f(x_{N,0}),f(x_{N,1}),\ldots,f(x_{N,N}))^{T}$, then $H_{N}v$ is an approximation to $(f'(x_{N,0}),f'(x_{N,1}),\ldots,f'(x_{N,N}))^{T}$. The second, components of the vector $w_{N}:=(w_{N,0},w_{N,1},\ldots,w_{N,N})^{T}\in\mathbb{R}^{N+1}$, are given by
\begin{equation*}
w_{N,j}:=\int_{0}^{l}\ell_{N,j}(x)\dd x,\quad j=0,1,\ldots,N,
\end{equation*}
and, for the same $v$ above, $w_{N}^{T}v$ is an approximation to $\int_{0}^{l}f(x)\dd x$. Both follow straightforwardly from approximating $f$ with the $N$-degree interpolating polynomial\begin{equation}\label{lagrange}
p_{N}(x):=\sum_{j=0}^{N}\ell_{N,j}(x)f(x_{N,j}),\quad x\in[0,l],
\end{equation}
which indeed satisfies $p(x_{N,i})=f(x_{N,i})$, $i=0,1,\ldots,N$. Even though at this moment it is not necessary to specify any particular choice of nodes, let us remark that in the case of Chebyshev extremal points, as we assume in Section \ref{s_results}, $H_{N}$ and $w_{N}$ can be obtained rather easily\footnote{Usually these quantities are provided with reference to the normalized interval $[-1,1]$.}. For these and other related aspects see, e.g., \cite{tref00}. See also \cite{dsv20} for a recent review of results on $H_{N}$.

\bigskip
The collocation method leads to a discrete version
\begin{equation}\label{dgeig}
B_{N}\Phi=\lambda M_{N}\Phi
\end{equation}
of the generalized eigenvalue problem \eqref{geig}, where the structure of the matrix representation of the finite-dimensional operators $B_{N}$ and $M_{N}$ depends on the specific model as detailed in the following Section \ref{s_discretizationA} for model A and Section \ref{s_discretizationB} for model B. In any case, the discrete problem has finite dimension, as it is posed on $X_{N}$.

Of course, as a main outcome one expects that the eigenvalues of \eqref{dgeig} approximate (part of) the eigenvalues of \eqref{geig}, the accuracy improving as $N$ increases. This is shown to be the case in Section \ref{s_results} by way of several numerical experiments, properly tuned to test the error behavior, to measure the convergence rate and to detect other peculiarities of the proposed approach. Note that the eigenvalues of \eqref{dgeig} can be computed with standard techniques for finite-dimensional generalized eigenvalue problems (e.g., Matlab's \texttt{eig}, based on the well-known QZ algorithm, see, e.g., \cite{gvl13}). Let us remark that we are anyway interested in the spectral radius of the next generation operator, so that we are mostly concerned with the dominant part of its spectrum.

Let us finally stress again that the use of the proposed methodology to approximate $R_{0}$ is based on the assumption that this number actually corresponds to a generalized eigenvalue. If this is not the case, already \eqref{geig} looses sense. Nevertheless, in Section \ref{s_results}, we report on some tests under the latter hypothesis, where the scheme is still able to give reasonable approximations (even though with slower convergence than what experimented in the case of generalized eigenvalues).
%-----------------------------------------------------------------------------
\subsection{Discretization of model A}
\label{s_discretizationA}
Let us recall the main ingredients \eqref{BA}, \eqref{MA} and \eqref{DMA} of the class of cell population models, model A in Section \ref{s_modelA}. The generalized eigenvalue problem \eqref{geig} reads
\begin{equation}\label{geigx}
(B\phi)(x)=\lambda(M\phi)(x),\quad x\in[0,l],
\end{equation}
with generalized eigenfunction $\phi\in\mathcal{D}(M)\setminus\{0\}$, i.e.,
\begin{equation}\label{BCA}
c(x)\phi(x)-D(x)\phi'(x)=0,\quad x=0,l.
\end{equation}
The proposed collocation consists in looking for an $N$-degree polynomial $p_{N}$ satisfying
\begin{equation*}
(Bp_{N})(x_{N,i})=\lambda(Mp_{N})(x_{N,i}),\quad i=1,\ldots,N-1,
\end{equation*}
together with
\begin{equation}\label{dBCA}
c(x_{N,i})p_{N}(x_{N,i})-D(x_{N,i})p_{N}'(x_{N,i})=0,\quad i=0,N.
\end{equation}
By using the Lagrange representation
\begin{equation}\label{plagrange}
p_{N}(x)=\sum_{j=0}^{N}\ell_{N,j}(x)\Phi_{j},\quad x\in[0,l],
\end{equation}
and by recalling the differentiation matrix $H_{N}$ introduced in Section \ref{s_discretization}, it is not difficult to recover \eqref{dgeig} with $B_{N}\in\mathbb{R}^{(N+1)\times(N+1)}$ given as
\begin{equation*}
B_{N;i,j}:=\left\{\setlength\arraycolsep{0.1em}\begin{array}{ll}
2\beta(x_{N,i}),&\quad i=j=1,\ldots,N-1,\\[2mm]
0,&\quad i=0,N \text{ or } i\neq j,
\end{array}
\right.
\end{equation*}
and $M_{N}\in\mathbb{R}^{(N+1)\times(N+1)}$ given as
\begin{equation*}
M_{N;i,j}:=\left\{\setlength\arraycolsep{0.1em}\begin{array}{lll}
C_{N;0,0}\delta_{0,j}-D_{N;0,0}h_{N;0,j},&\quad i=0,&\;j=0,1,\ldots,N,\\[2mm]
[H_{N}(C_{N}-D_{N}H_{N})+\Sigma_{N}]_{i,j},&\quad i=1,\ldots,N,&\;j=0,1,\ldots,N,\\[2mm]
C_{N;N,N}\delta_{N,j}-D_{N;N,N}h_{N;N,j},&\quad i=N,&\;j=0,1,\ldots,N,
\end{array}\right.
\end{equation*}
where $\delta_{i,j}$ is the Kronecker's delta and $C_{N},D_{N},\Sigma_{N}\in\mathbb{R}^{(N+1)\times(N+1)}$ are defined as
\begin{equation*}
C_{N;i,j}:=\left\{\setlength\arraycolsep{0.1em}\begin{array}{ll}
c(x_{N,i}),&\quad i=j=0,1,\ldots,N,\\[2mm]
0,&\quad i\neq j,
\end{array}
\right.
\end{equation*}
\begin{equation*}
D_{N;i,j}:=\left\{\setlength\arraycolsep{0.1em}\begin{array}{ll}
D(x_{N,i}),&\quad i=j=0,1,\ldots,N,\\[2mm]
0,&\quad i\neq j,
\end{array}
\right.
\end{equation*}
\begin{equation*}
\Sigma_{N;i,j}:=\left\{\setlength\arraycolsep{0.1em}\begin{array}{ll}
\beta(x_{N,i})+\mu(x_{N,i}),&\quad i=j=0,1,\ldots,N,\\[2mm]
0,&\quad i\neq j.
\end{array}
\right.
\end{equation*}
Note how the first and last rows of $B_{N}$ and $M_{N}$ realize the discrete boundary conditions \eqref{dBCA}, simulating those characterizing $\mathcal{D}(M)$ in \eqref{DMA}, i.e., \eqref{BCA}.
%-----------------------------------------------------------------------------
\subsection{Discretization of model B}
\label{s_discretizationB}
Let us recall the main ingredients \eqref{BB}, \eqref{MB} and \eqref{DMB} of the class of age-structured epidemics, model B in Section \ref{s_modelB}. The generalized eigenvalue problem \eqref{geig} reads again as \eqref{geigx}, but with generalized eigenfunction $\phi\in\mathcal{D}(M)\setminus\{0\}$, i.e., satisfying
\begin{equation}\label{BCB}
\phi(0)=\theta\int_{0}^{l}\beta(x)\phi(x)\dd x.
\end{equation}
The proposed collocation consists in looking for an $N$-degree polynomial $p_{N}$ satisfying
\begin{equation*}
(Bp_{N})(x_{N,i})=\lambda(Mp_{N})(x_{N,i}),\quad i=1,\ldots,N,
\end{equation*}
together with
\begin{equation}\label{dBCB}
p_{N}(0)=\theta\int_{0}^{l}\beta(x)p_{N}(x)\dd x,
\end{equation}
approximating the integral above as well as the one defining $B$ in \eqref{BB} by quadrature. By using again the Lagrange representation \eqref{plagrange} and by recalling the (vector of) quadrature weights $w_{N}$ introduced in Section \ref{s_discretization}, it is not difficult to recover \eqref{dgeig} with $B_{N}\in\mathbb{R}^{(N+1)\times(N+1)}$ given as
\begin{equation*}
B_{N;i,j}:=\left\{\setlength\arraycolsep{0.1em}\begin{array}{lll}
0,&\quad i=0,&\;j=0,1,\ldots,N,\\[2mm]
w_{N,j}K(x_{N,i},x_{N,j}),&\quad i=1,\ldots,N,&\;j=0,1,\ldots,N,
\end{array}
\right.
\end{equation*}
and $M_{N}\in\mathbb{R}^{(N+1)\times(N+1)}$ given as
\begin{equation*}
M_{N;i,j}:=\left\{\setlength\arraycolsep{0.1em}\begin{array}{lll}
\delta_{0,j}-\theta w_{N,j}\beta(x_{N,j}),&\quad i=0,&\;j=0,1,\ldots,N,\\[2mm]
h_{N;i,j}+[\gamma(x_{N,i})+\delta(x_{N,i})]\delta_{i,j},&\quad i=1,\ldots,N,&\;j=0,1,\ldots,N,
\end{array}
\right.
\end{equation*}
where $\delta_{i,j}$ is again the Kronecker's delta. Note how the first row of both $B_{N}$ and $M_{N}$ realizes the (quadrature of the) discrete boundary condition \eqref{dBCB} simulating the one characterizing $\mathcal{D}(M)$ in \eqref{DMB}, i.e., \eqref{BCB}.
%-----------------------------------------------------------------------------
%-----------------------------------------------------------------------------
\section{Results}
\label{s_results}
A first series of numerical tests is presented in Section \ref{s_results1} with the aim of illustrating the convergence properties of the proposed techniques, as well as related aspects such as, e.g., the effect of the smoothness of the models' coefficients. Then in Section \ref{s_results2} we use the tested algorithms to perform some quantitative studies in the context of varying parameters, checking whether the outcomes confirm the theoretical expectations.

All the experiments are run on a MacBook Pro 2.3 GHz Intel Core i7 with 16 GB memory, through codes written in Matlab version R2019a (freely available at \url{http://cdlab.uniud.it/software} as anticipated in the Introduction). In these codes we always use Chebyshev extremal points as collocation nodes, i.e., $x_{N,i}=l[1-\cos(i\pi/N)]/2$, $i=0,1,\ldots,N$. We refer to \cite{tref00} for their numerous advantageous properties in the context of numerical interpolation, differentiation and quadrature.
%-----------------------------------------------------------------------------
\subsection{Numerical properties}
\label{s_results1}
The numerical properties of the proposed approach are mainly investigated by analyzing the behavior of the error $|R_{0,N}-R_{0}|$ between the approximated and the exact basic reproduction numbers. Here for {\it exact} we mean either the theoretical value of $R_{0}$ when this is explicitly available, or a reference value $R_{0,\bar N}$ computed with a suitably large $\bar N$ when an exact value is not attainable.

A second error is also presented when possible, viz. $\|p_{N}-\phi\|_{\infty,M}$. This is the error between the exact generalized eigenfunction (if available) and its collocation approximation, measured as the maximum absolute value of their differences on a mesh of $M$ equidistant points in $[0,l]$. In particular, we always use $M=1\,000$. The collocation polynomial $p_{N}$ is reconstructed from the computed generalized eigenvector $\Phi$ associated to the dominant generalized eigenvalue of \eqref{dgeig} through barycentric interpolation by applying the algorithm proposed in \cite{bertre04}. Of course, the same normalization is prescribed from time to time for both $p_{N}$ and $\phi$.

Finally, the discrete generalized eigenvalue problem \eqref{dgeig} is solved either by Matlab's \texttt{eig} or \texttt{eigs}. We also test the discretization
\begin{equation}\label{deig}
B_{N}M_{N}^{-1}\Psi=\lambda\Psi,
\end{equation}
of the eigenvalue problem
\begin{equation}\label{eig}
BM^{-1}\psi=\lambda\psi
\end{equation}
equivalent to \eqref{geig}, where $\psi=M\phi$ for $\phi$ in \eqref{geig} and $\Psi=M_{N}\Phi$ for $\Phi$ in \eqref{dgeig}. These eigenproblems are implemented as \texttt{eig(B,M)} for \eqref{dgeig} and \texttt{eig(B/M)} for \eqref{deig} (or the alternative versions through \texttt{eigs}).
%-----------------------------------------------------------------------------
\subsubsection{Specific models}
\label{s_results1models}
As far as the models to test are concerned, below we list all the specific choices by giving the defining rates and coefficients. In parallel we also give some of their key features.

For all the instances of model A we set $l=1$, $c(x)=l-x$, $x\in[0,l]$, and
\begin{equation*}
D(x):=\tilde D\cdot\left[\frac{4}{l^{2}}x(l-x)+1\right]\geq\tilde D,\quad x\in[0,l].
\end{equation*}
The other ingredients differ as described next.
\begin{enumerate}[label=(A\arabic*),ref=(A\arabic*)]
\item\label{A1} [the {\it immortal} case] Let $\tilde D=1$, $\mu\equiv0$ and
\begin{equation*}
\beta(x):=\tilde\beta\cdot\left[\frac{27}{2l^{3}}x^{2}(l-x)+1\right],\quad x\in[0,l],
\end{equation*}
for $\tilde\beta=1$. With these choices the next generation operator is compact according to \cite{brv20} and thus $R_{0}$ is a generalized eigenvalue. It is not difficult to recover $R_{0}=2$ from \eqref{upperR0eig}. Moreover, some calculations that we omit allow to show, starting from \eqref{geig}, that the corresponding generalized eigenfunction is
\begin{equation}\label{phiA1}
\phi(x)=e^{\int_{0}^{x}\frac{c(y)}{D(y)}\dd y},\quad x\in[0,l],
\end{equation}
normalized as $\phi(0)=1$, and the eigenfunction in \eqref{eig} is $\psi=\beta\phi$.
\item\label{A2} [the {\it proportional} case] Let $\tilde D=1$,
\begin{equation*}
\mu(x):=\tilde\mu\frac{x}{l},\quad \beta(x):=\tilde\beta\mu(x),\quad x\in[0,l],
\end{equation*}
for $\tilde\mu=1$ and $\tilde\beta=1.5$. Again, it turns out that we are in the compact case and that $R_{0}=2\tilde\beta/(\tilde\beta+1)=1.2$ follows from \eqref{upperR0eig}. Moreover, the corresponding generalized eigenfunction is again given by \eqref{phiA1}, and the eigenfunction in \eqref{eig} is $\psi=(\tilde{\beta}+1)\mu\phi$.
\item\label{A3} [the {\it general} case] Let $\mu$ be the same as in case \ref{A2} and $\beta$ be expressed as in case \ref{A1}, but with $\tilde\beta=10$. We consider either
\begin{enumerate}[label=(A3.\arabic*),ref=(A3.\arabic*)]
\item\label{A31} the compact case: $\tilde D=1$;   \footnote{Here we refer to the compact case  since when there is diffusion everywhere in space and under technical assumptions on the model parameters, the next generation operator is compact according to \cite{brv20}.}
\item\label{A32} a ``almost non-compact'' case: $\tilde D=10^{-6}$;
\item\label{A33} the non-compact case: $\tilde D=0$.  \footnote{For the case of lack of diffusion $R_{0}$ can be computed analytically and it turns out to be not an eigenvalue.}
\end{enumerate}
Independently of these choices, $R_{0}$ is unknown, as well as the corresponding generalized eigenfunction.
\end{enumerate}

For all the instances of model B we set again $l=1$ and the rest is defined below.
\begin{enumerate}[label=(B\arabic*),ref=(B\arabic*)]
\item\label{B1} [the {\it age-independent} case] Let $\mu\equiv\tilde\mu=28$ (so that $e^{-\tilde\mu l}\simeq0$), $\beta_{0}\equiv\mu$ (so that \eqref{intbeta0Pi=1} holds), $\gamma\equiv\tilde\gamma=1$ and $\delta\equiv\tilde\delta=1$ and
\begin{equation*}
K(x,y):=\tilde k\tilde\mu e^{-\tilde\mu y},\quad x,y\in[0,l],
\end{equation*}
for $\tilde k=52$. It also follows that
\begin{equation*}
\beta(x)=\tilde\mu e^{-\tilde\mu x},\quad \Pi_{1}(x)=e^{-(\tilde\gamma+\tilde\delta)x},\quad x\in[0,l].
\end{equation*}
Finally, let $\theta=1/7$. With these choices the next generation operator is compact according to \cite{brv20}. Due to the finite maximum age there is a very small fraction of immortal individuals, but still we can explicitly obtain $R_{0}\simeq\tilde{k}/[(1-\theta)\tilde\mu+\tilde\gamma+\tilde\delta]=2$. It is evident from \eqref{NG-Age-Epidemic} that the corresponding eigenfunctions $\psi$ in \eqref{eig} are the constant functions given that $K$ is independent of $x$. Then one can recover through \eqref{MB} the generalized eigenfunctions
\begin{equation}\label{phiB1}
\phi(x)=e^{-2x}\left(\phi(0)-\frac{\psi}{2}\right)+\frac{\psi}{2},\quad x\in[0,l],
\end{equation}
with $\phi(0)\simeq\psi/182$ following from \eqref{BCB}, with $\psi$ playing the role of normalizing factor.
\item\label{B2} [the {\it proportionate} case] Let
\begin{equation*}
K(x,y):=\tilde k x^{2}(l-x)^{2}\cdot\frac{\Pi_{0}(y)}{\int_{0}^{l}\Pi_{0}(z)\dd z},\quad x,y\in[0,l],
\end{equation*}
for $\tilde k=16\,065/64$ and
\begin{equation}\label{Pi0B2}
\Pi_{0}(x):=\left(\frac{l-x}{l}\right)^{\alpha},\quad x\in[0,l],
\end{equation}
for some $\alpha>0$. Note that $\int_{0}^{l}\Pi_{0}(z)\dd z=\frac{l}{\alpha+1}$. Let moreover $\gamma$ and $\delta$ satisfy
\begin{equation}\label{gammadelta}
\gamma(x)+\delta(x)=\frac{1}{l-x},\quad x\in[0,l],\footnote{For this specific choice, and for the sake of pragmatism, here we exclude the last node $x_{N,N}=l$ from the computations to avoid overflow. In \cite{brv20} the issue is tackled more rigorously from the numerical standpoint.}
\end{equation}
so that also $\Pi_{1}(x)=\frac{l-x}{l}$, $x\in[0,l]$, follows. We assume absence of vertical transmission, i.e., $\theta=0$, which makes useless to specify $\beta$. With these choices the next generation operator is again compact according to \cite{brv20}. Since $K$ above is the product of functions in each of the variables $x$ and $y$, the next generation operator in \eqref{NG-Age-Epidemic} becomes a rank one operator and we can explicitly compute
\begin{equation*}
R_{0}=\frac{2\tilde k(\alpha+1)l^{5}}{(\alpha+2)(\alpha+4)(\alpha+5)(\alpha+6)}.
\end{equation*}
As far as $\alpha$ is concerned, we consider
\begin{enumerate}[label=(B2.\arabic*),ref=(B2.\arabic*)]
\item\label{B21} $\alpha=1/4$: life expectancy is $0.8l$ and $R_{0}=2$;
\item\label{B22} $\alpha=1$: life expectancy is $0.5l$ and $R_{0}=1.59375$.
\end{enumerate}
Finally, and independently of $\alpha$, the eigenfunction corresponding to $R_{0}$ is $\psi(x)=x^{2}(l-x)^{2}$, normalized as $\psi(l/2)=l^{4}/16$, and the generalized eigenfunction is
\begin{equation*}
\phi(x)=\Pi_{1}(x)\int_{0}^{x}\frac{y^{2}(l-y)^{2}}{\Pi_{1}(y)}\dd y,\quad x\in[0,l],
\end{equation*}
which, for the choices above, turns out to be a polynomial of degree $5$, viz.
\begin{equation}\label{phiB2}
\phi(x)=\frac{x^{3}(l-x)(4l-3x)}{12},\quad x\in[0,l]
\end{equation}
(normalized as $\phi(l/2)=5l^{5}/384$).
\item\label{B3} [the {\it general} case] Let
\begin{equation*}
K(x,y):=\tilde k e^{-\left(\frac{x-y}{l_{0}}\right)^{2}}\cdot\frac{\Pi_{0}(y)}{\int_{0}^{l}\Pi_{0}(z)\dd z},\quad x,y\in[0,l],
\end{equation*}
for $\tilde k=25$, $l_{0}=0.1l$ and $\Pi_{0}(x):=e^{-\frac{\alpha x}{l-x}}$, $x\in[0,l]$, for $\alpha=0.1$ (life expectancy is $0.8l$). Some calculations give
\begin{equation*}
\int_{0}^{l}\Pi_{0}(x)\dd x=l\left(1-\alpha e^{\alpha}\int_{\alpha}^{+\infty}\frac{e^{-x}}{x}\dd x\right)
\end{equation*}
Let moreover $\gamma$ and $\delta$ be as in case \ref{B1}. Here, infection events occur mostly between individuals of the same age and there are chronic infected individuals since $e^{-(\tilde\gamma+\tilde\delta)l}>0$. Finally, let $\theta=0.5$ and $\beta(x):=b(x)\Pi_{0}(x)/\int_{0}^{l}b(y)\Pi_{0}(y)\dd y$, $x\in[0,l]$, for $b(x):=(x/l)^{2}e^{-6x/l}$, which ensures \eqref{intbeta0Pi=1} and maximum fertility at $x=l/3$. For these choices the next generation operator is compact \footnote{In general, under technical assumptions on the infection kernel, the next generation operator turns out to be compact according to \cite{brv20}.}, but both $R_{0}$ and the corresponding generalized eigenfunction are unknown.
\end{enumerate}
%-----------------------------------------------------------------------------
\subsubsection{Results for model A}
\label{s_results1A}
\begin{figure}[!ht]
\centering
\includegraphics[width=.48\textwidth]{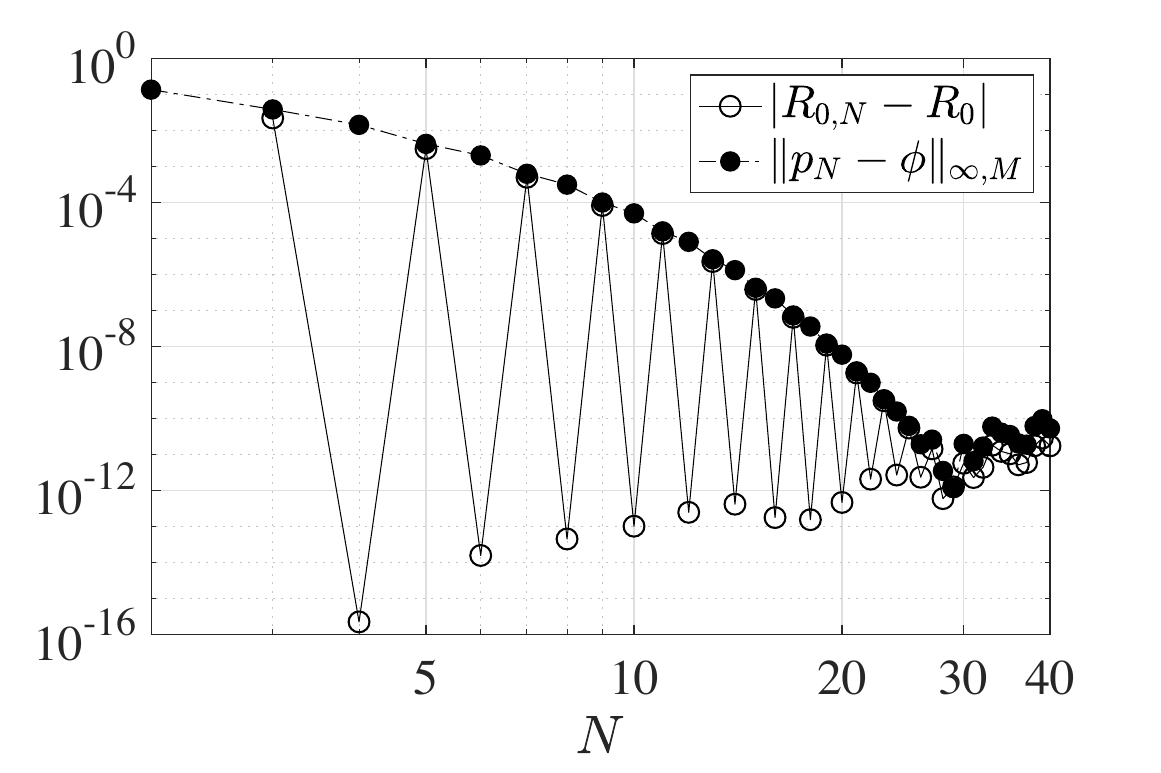}
\includegraphics[width=.48\textwidth]{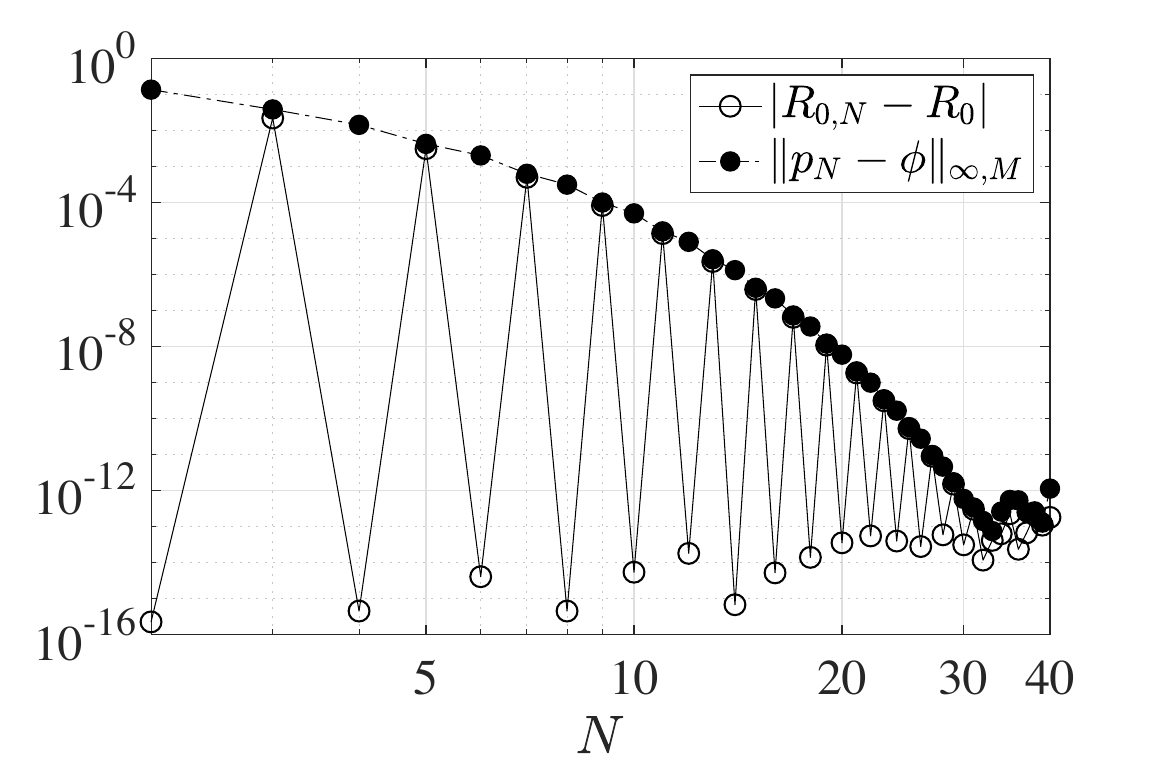}
\caption{case \ref{A1} -- errors for increasing $N$ on $R_{0}=2$ and on $\phi$ in \eqref{phiA1}, computed as \eqref{geig} through \texttt{eig(B,M)} (left) and as \eqref{eig} through \texttt{eig(B/M)} (right).}
\label{f_immortal_old}
\end{figure}

\begin{figure}[!ht]
\centering
\includegraphics[width=.48\textwidth]{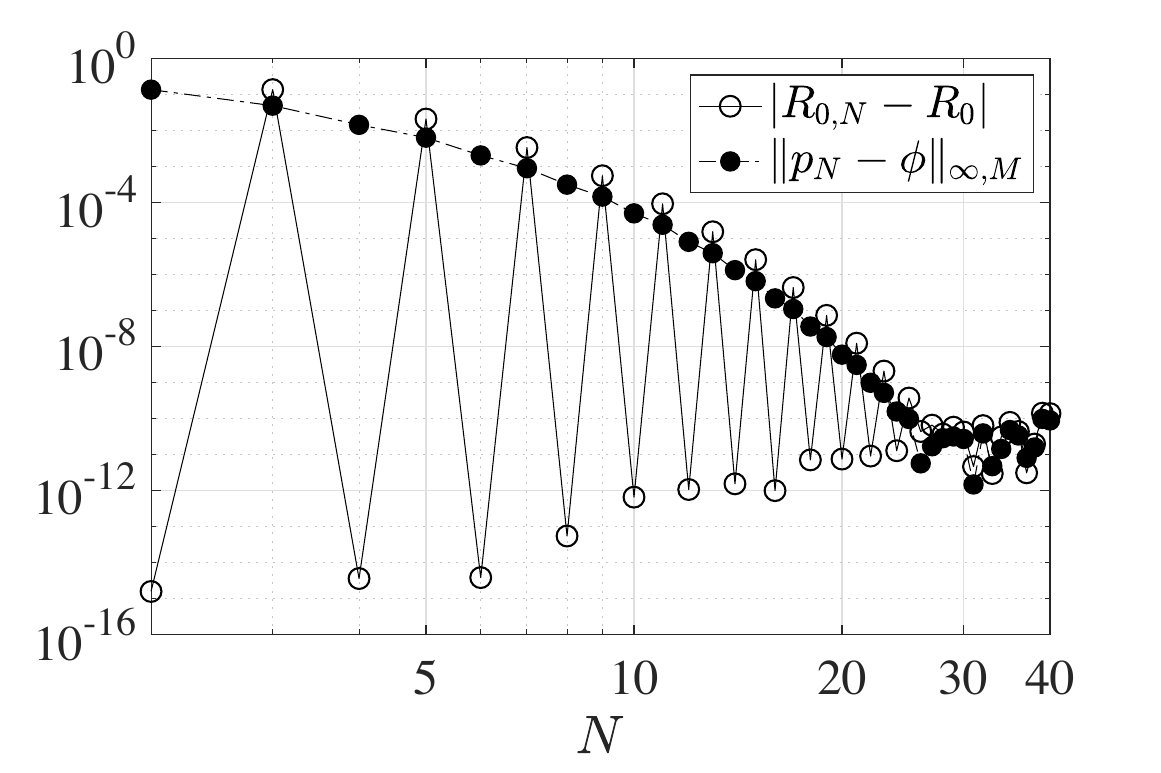}
\includegraphics[width=.48\textwidth]{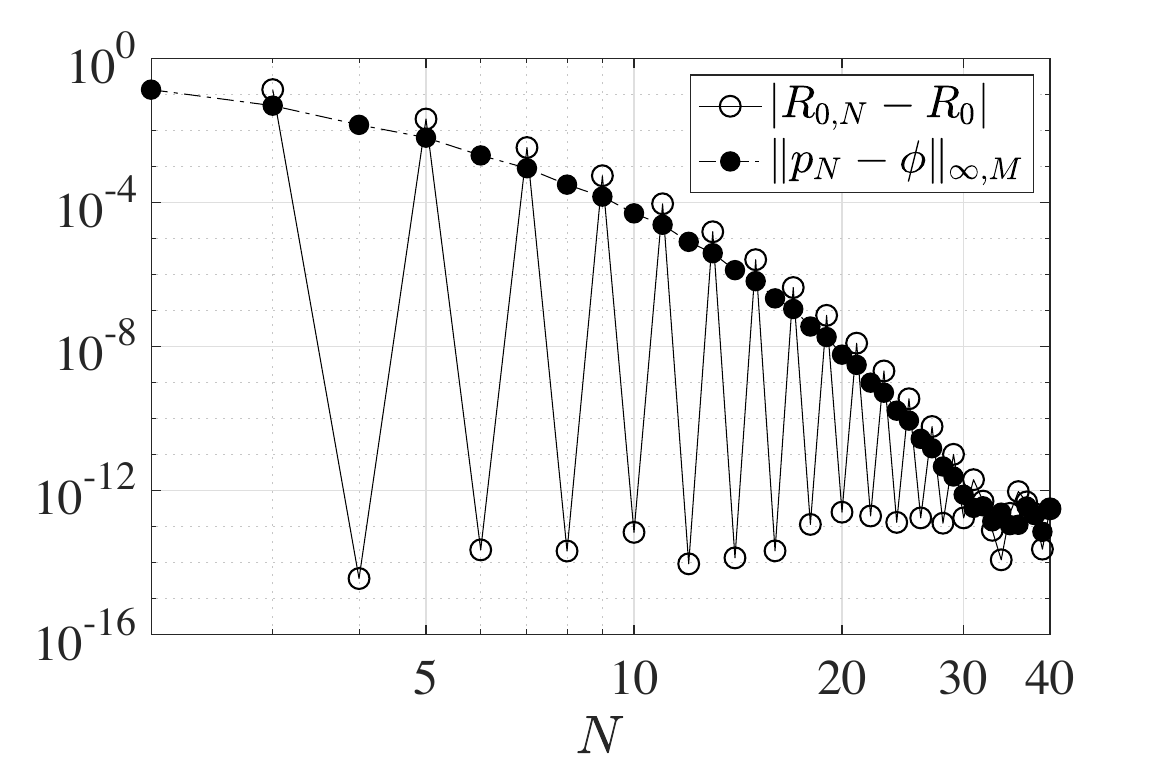}
\caption{case \ref{A2} -- errors for increasing $N$ on $R_{0}=1.2$ and on $\phi$ in \eqref{phiA1}, computed as \eqref{geig} through \texttt{eig(B,M)} (left) and as \eqref{eig} through \texttt{eig(B/M)} (right).}
\label{f_proportional_old}
\end{figure}

As a first experiment we test the convergence properties of the proposed approach for case \ref{A1}. The spectrally accurate behavior (namely the error decays faster than $O(N^{-k})$ for any natural $k$, \cite{tref00}) is evident in Figure \ref{f_immortal_old} (left) for the approximation of both $R_{0}$ and the relevant generalized eigenfunction $\phi$. To note that for the former the approximation is already very good for low even values of $N$, due to the (anti-) symmetry properties of the Chebyshev differentiation matrix $H_{N}$ \cite{tref00}. Similar trends emerge also in the right panel, where \eqref{eig} is solved through \texttt{eig(B/M)} instead of solving \eqref{geig} through \texttt{eig(B,M)} as in the left panel. Let us just remark that the choice of \texttt{eig(B/M)} seems slightly better from the point of view of algorithmic stability. Finally, let us mention that we have compared also with the use of \texttt{eigs}, obtaining indistinguishable results. As these outcomes are unchanged in the remaining experiments, we show only the results relevant to the use of \texttt{eig(B/M)}, although we measure always the error on the generalized eigenfunction $\phi$.
Figure \ref{f_proportional_old} concerns the same analysis above, but for case \ref{A2}: the results are qualitatively the same.

\begin{figure}[!ht]
\centering
\includegraphics[width=.48\textwidth]{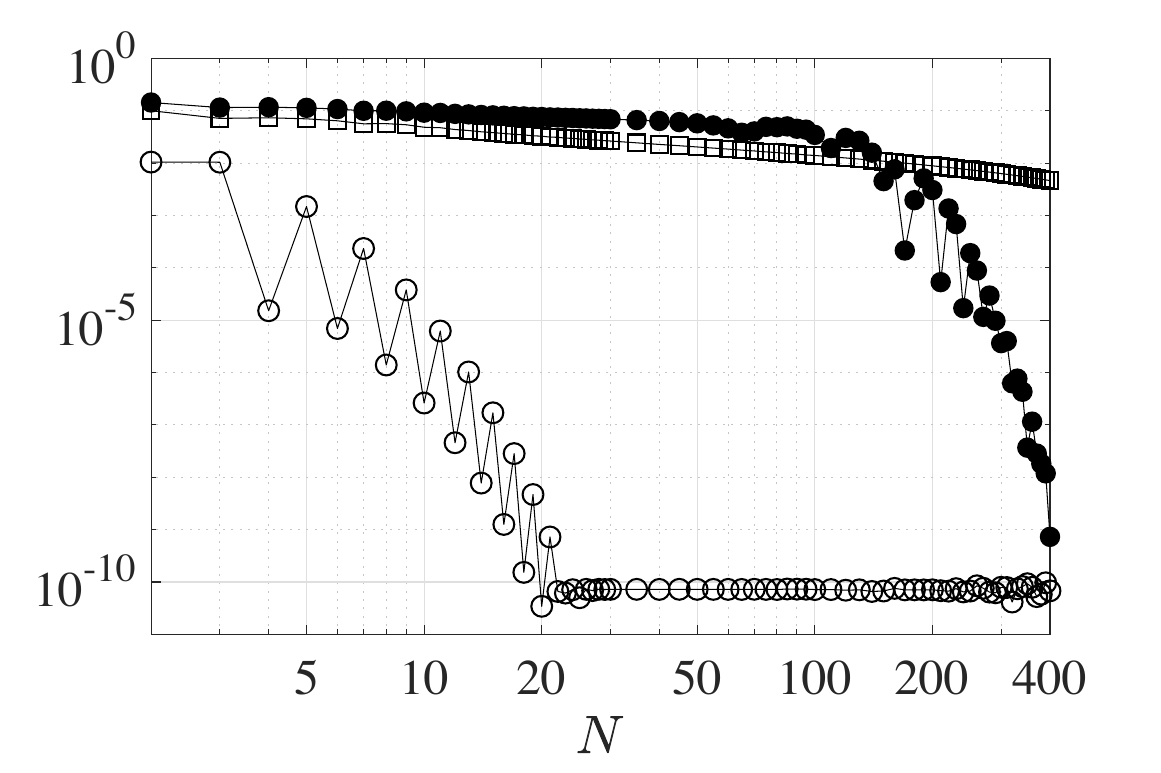}
\caption{case \ref{A3} -- error for increasing $N$ on $R_{0}\simeq1.953037627191416$ ($\circ$, case \ref{A31}), $R_{0}\simeq1.820082420868451$ ($\bullet$, case \ref{A32}), $R_{0}\simeq1.863455805786607$ ($\square$, case \ref{A33}) computed as \eqref{eig} through \texttt{eig(B/M)} (reference values $R_{0,\bar N}$ computed with $\bar N=1\,000$).}
\label{f_generalA}
\end{figure}

In the previous tests both $R_{0}$ and the relevant generalized eigenfunction are analytically known. As the same is not possible for case \ref{A3}, we show in Figure \ref{f_generalA} only the error on $R_{0}$ with respect to a reference value $R_{0,N}$ as explained at the beginning of Section \ref{s_results1}. As expected, the convergence is spectrally accurate for case \ref{A31}, in which $\tilde D=1$ ensures compactness of the next generation operator (according to \cite{brv20}). Theoretically, also $\tilde D=10^{-6}\neq0$ in case \ref{A32} guarantees compactness, but it is clearly visible from the plot that much larger values of $N$ are necessary to start appreciating the spectral accuracy. It is indeed reasonable to expect that the value of $\tilde D$ affects the error constants, causing their increase as $\tilde D\rightarrow0$, still being the problem compact. When we deal instead with case \ref{A33}, in which the absence of diffusion causes the loss of compactness, convergence still occurs, even though at a fixed rate (seemingly linear). The fact is somehow surprising (and certainly merits future investigation), given that without compactness \eqref{geig} may even become meaningless and we are thus using a finite-dimensional eigenvalue problem to approximate components of the spectrum possibly other than the point one.

As a final remark for model A, let us mention that all the eigenfunctions computed with the proposed approach in case \ref{A3} turned out to approximate seemingly smooth curves (as they are expected to be). To avoid redundancy, we give examples only for case \ref{B3} in Section \ref{s_results2}.
%-----------------------------------------------------------------------------
\subsubsection{Results for model B}
\begin{figure}[!ht]
\centering
\includegraphics[width=.48\textwidth]{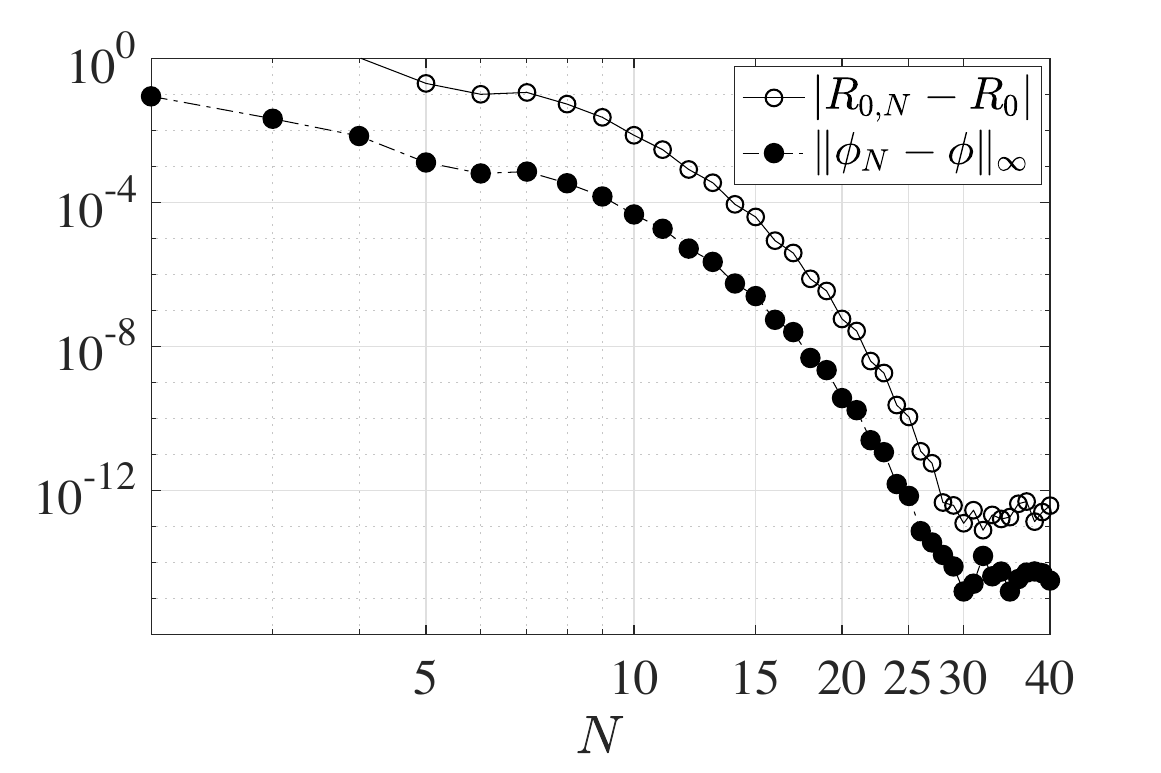}
\caption{case \ref{B1} -- errors for increasing $N$ on $R_{0}\simeq2$ and on $\phi$ in \eqref{phiB1}, computed as \eqref{eig} through \texttt{eig(B/M)} (reference value $R_{0,\bar N}=1.999999999982096$ computed with $\bar N=1\,000$).}
\label{f_ageindep}
\end{figure}

As far as case \ref{B1} is concerned, the results of Figure \ref{f_ageindep} show spectral accuracy in the approximation of both $R_{0}$ and the relevant generalized eigenfunction as expected.
\label{s_results1B}

Slightly dissimilar is the situation for case \ref{B2}, where we still get convergence, but with trends different from what experimented so far and, moreover, possibly depending on $\alpha$. First of all, let us recall that the generalized eigenfunction \eqref{phiB2} is a polynomial of degree $5$, justifying the sudden drop of the error to machine precision occurring with $N=6$ in both the left and right panels. The theoretical proof of convergence fully elaborated in \cite{brv20} shows that this indeed happens at $N=6$ and not already at $N=5$ as one may expect. As far as the approximation of $R_{0}$ is concerned, instead, the same behavior is ensured if enough regularity is provided, as it is the case for $\alpha=1$ in the right panel. In fact, for instance, it is evident from \eqref{Pi0B2} that $\Pi_{0}$ is smooth for $\alpha\geq1$, but for $\alpha<1$ it becomes rational and blows up at $x=l$. Indeed, the choice $\alpha=1/4$ in the left panel prevents the method to perform the standard spectral accuracy, and convergence of fixed order (seemingly $4$) occurs. That the regularity of the model ingredients plays a role in the convergence analysis is a general fact, but for the specific problem currently investigated we leave the elaboration of the necessary details to \cite{brv20}.
\begin{figure}[!ht]
\centering
\includegraphics[width=.48\textwidth]{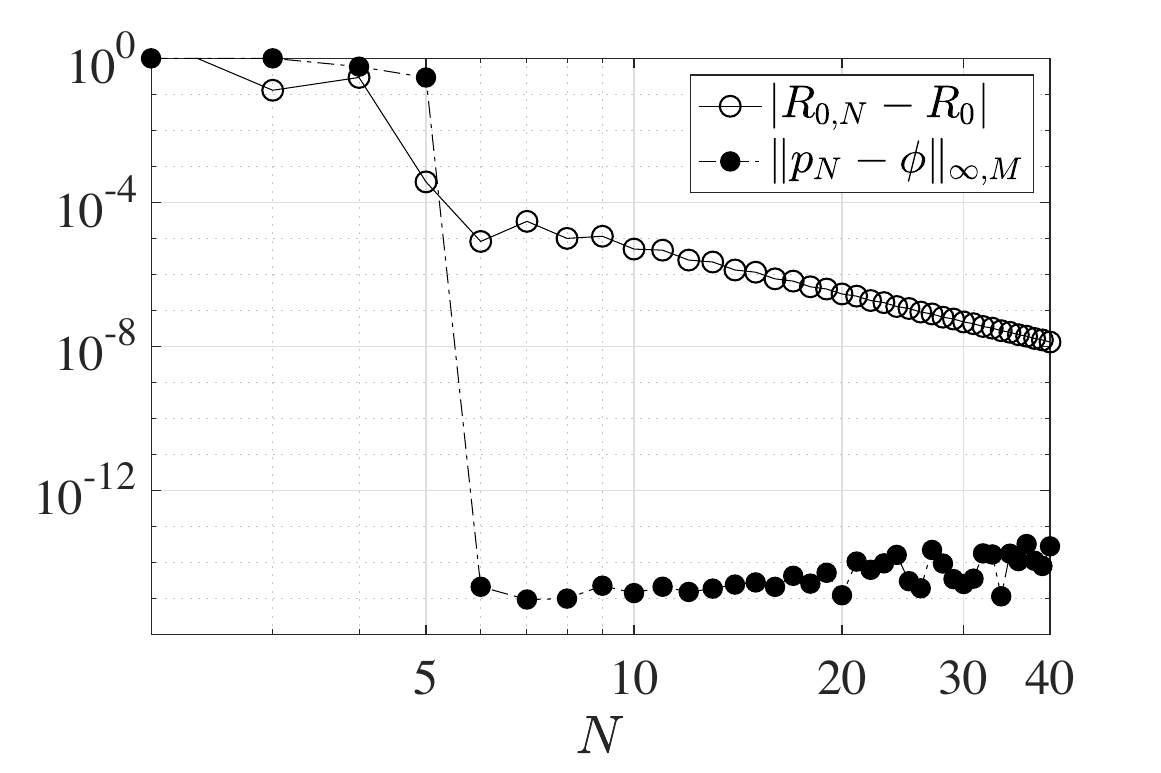}
\includegraphics[width=.48\textwidth]{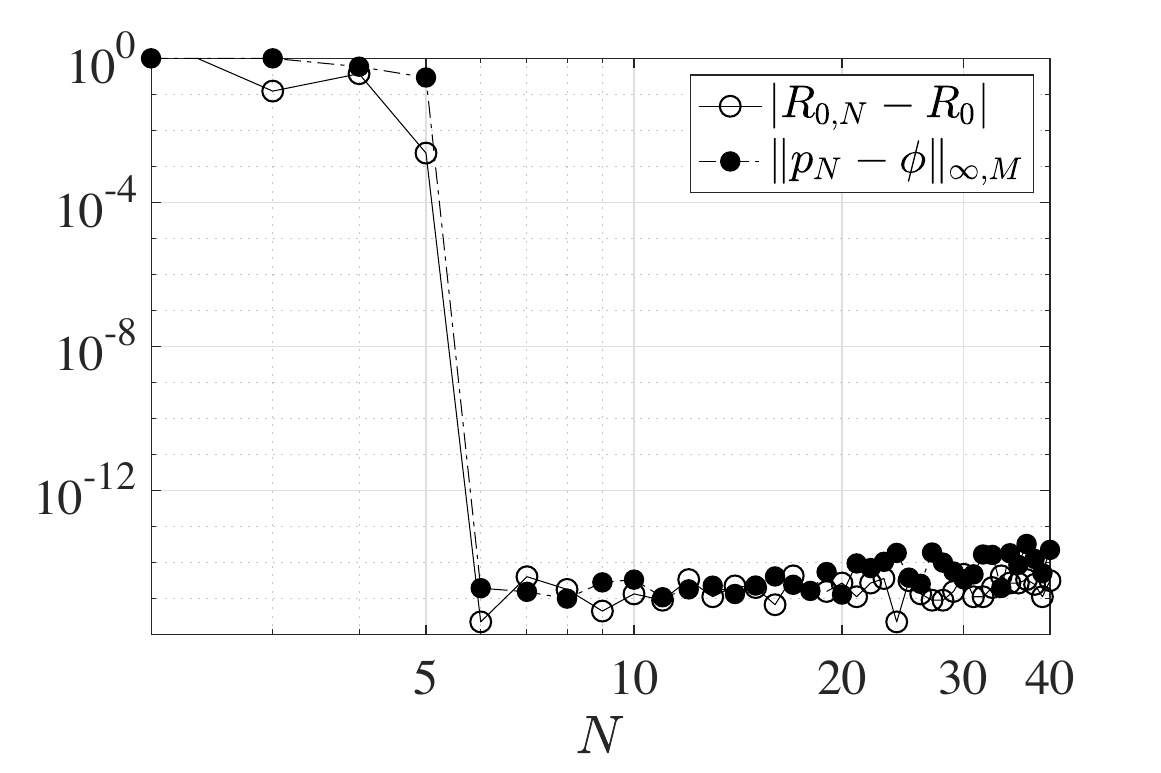}
\caption{case \ref{B2} -- errors for increasing $N$ on $R_{0}$ and on $\phi$ in \eqref{phiB2}, computed as \eqref{eig} through \texttt{eig(B/M)} for $\alpha=1/4$ (left, case \ref{B21}, $R_{0}=2$) and $\alpha=1$ (right, case \ref{B22}, $R_{0}=1.59375$).}
\label{f_proportionate}
\end{figure}

\begin{figure}[!ht]
\centering
\includegraphics[width=.48\textwidth]{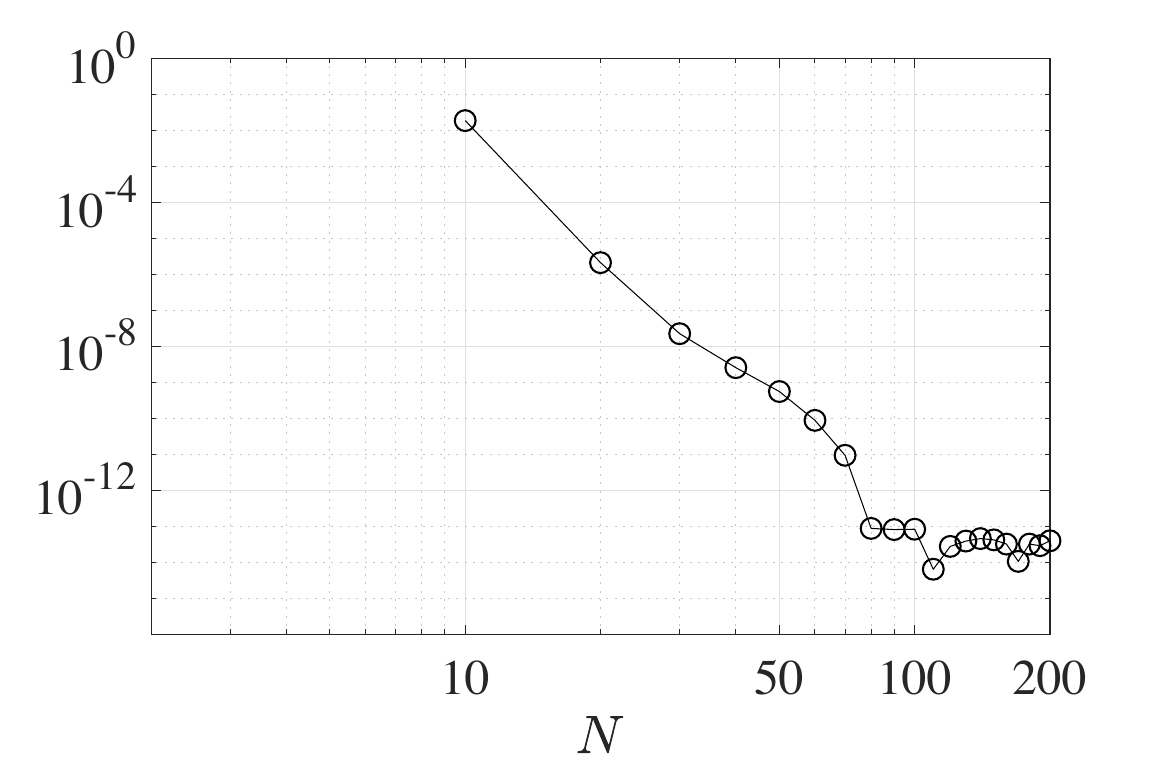}
\caption{case \ref{B3} -- error for increasing $N$ on $R_{0}$, computed as \eqref{eig} through \texttt{eig(B/M)} (reference value $R_{0,\bar N}=1.403437609535222$ computed with $\bar N=1\,000$).}
\label{f_generalB}
\end{figure}

Finally, the analysis of case \ref{B3} reveals spectral accuracy again, Figure \ref{f_generalB}. Recall that in this case, as for case \ref{A3} of the preceding section, neither $R_{0}$ nor the relevant generalized eigenfunction are known.
%-----------------------------------------------------------------------------
\subsection{Robust analysis}
\label{s_results2}
To conclude, we give some examples of practical application of the proposed techniques to analyze the behavior of the chosen models in the presence of varying parameters. Below, the choices of the discretization index $N$ are instructed by the results discussed in the previous section. Moreover, the intervals of variation of the concerned parameters are discretized by using uniform meshes of $P$ points. Corresponding values of $R_{0}$ are repeatedly approximated with the chosen $N$ for each of these points, or of their Cartesian product in the presence of two varying parameters. In the latter case we adopt standard contouring algorithms, viz. Matlab's \texttt{contourf}. Besides the graphical output, we give also some indication of the overall computational time with respect to the parameter(s) grid size $P$. Let us anticipate that all the tests concern the general instances of either model A or B, namely cases \ref{A3} or \ref{B3} as described in Section \ref{s_results1}.

\begin{figure}[!ht]
\centering
\includegraphics[width=.48\textwidth]{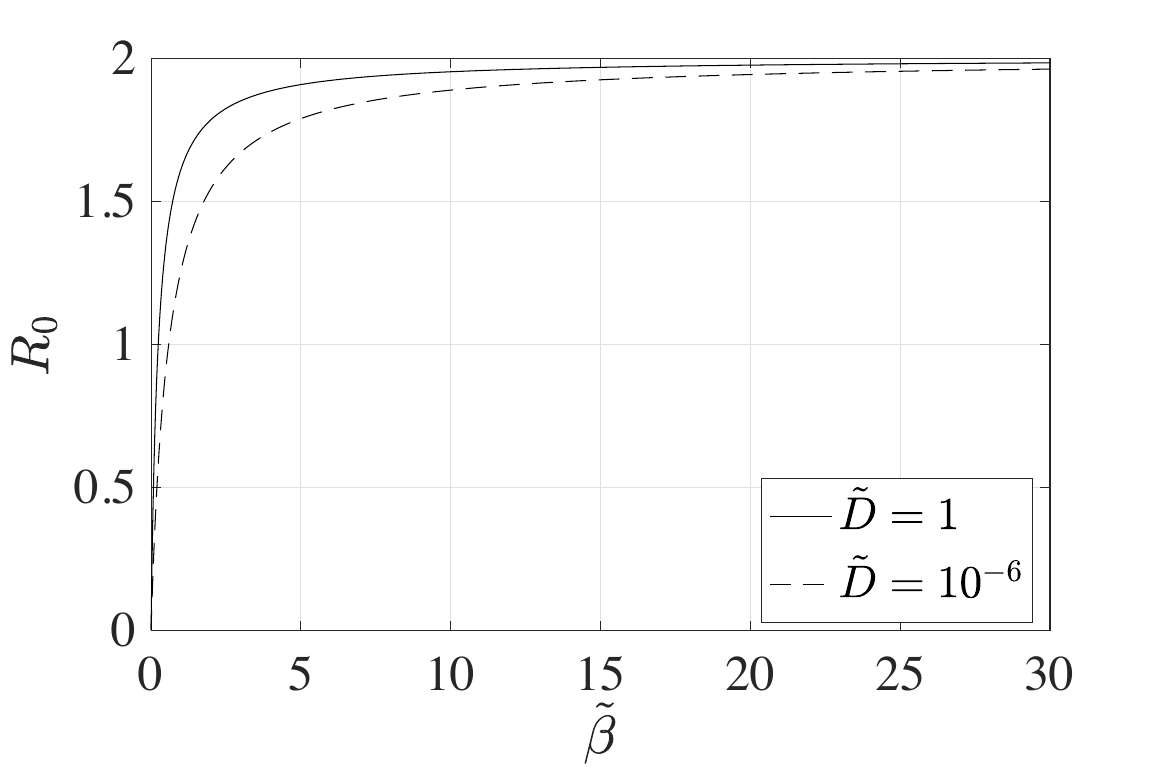}
\includegraphics[width=.48\textwidth]{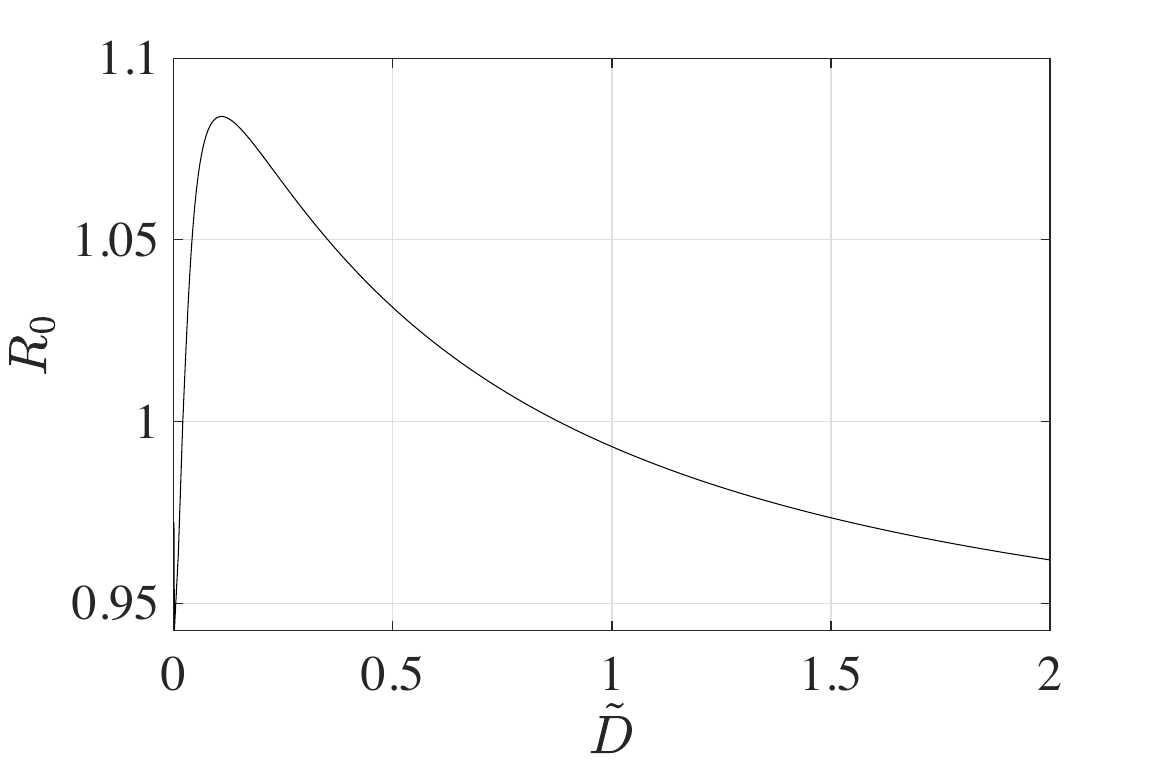}
\caption{case \ref{A3} -- $R_{0}$ as a function of $\tilde\beta$ for varying $\tilde D$ and $\tilde\mu=1$ (left, values $R_{0,N}$ computed with $N=30$ for case \ref{A31} and $N=500$ for case \ref{A32}) and $R_{0}$ as a function of $\tilde D$ for $\tilde\beta=10$ and $\tilde\mu=50$ (right, values $R_{0,N}$ computed with $N=300$).}
\label{f_vs_R0}
\end{figure}

In Figure \ref{f_vs_R0} we investigate for case \ref{A3} how $R_{0}$ varies as a function of either global fertility parameter $\tilde\beta$ (left panel) or global diffusion parameter$\tilde D$ (right panel). For the former $\tilde\mu=1$ and both cases \ref{A31} and \ref{A32} are considered. For the latter $\tilde\mu=50$ and $\tilde\beta=10$. As expected, $R_{0}$ grows with respect to the fertility from zero up to the asymptote $R_{0}=2$, since we assumed symmetric division \eqref{R0<2}, with larger diffusion values favoring the speed of increase. Instead interestingly, at fixed fertility, there is a single peak with respect to the diffusion, i.e.  there exist an intermediate diffusion maximizing the basic reproduction number. At this point, a naive approach is that a population adopting such strategy would persist and become an unbeatable population.\footnote{The role of the diffusion from the evolutionary point of view is out of the scope of the present work.}
Computationally speaking, both figures are obtained by repeatedly computing $R_{0}$ with respect to the varying parameter upon discretizing the latter on a mesh of $P=1\,000$ points, procedure that took on average just few seconds (say $1$ to $5$).

\begin{figure}[!ht]
\centering
\includegraphics[width=.48\textwidth]{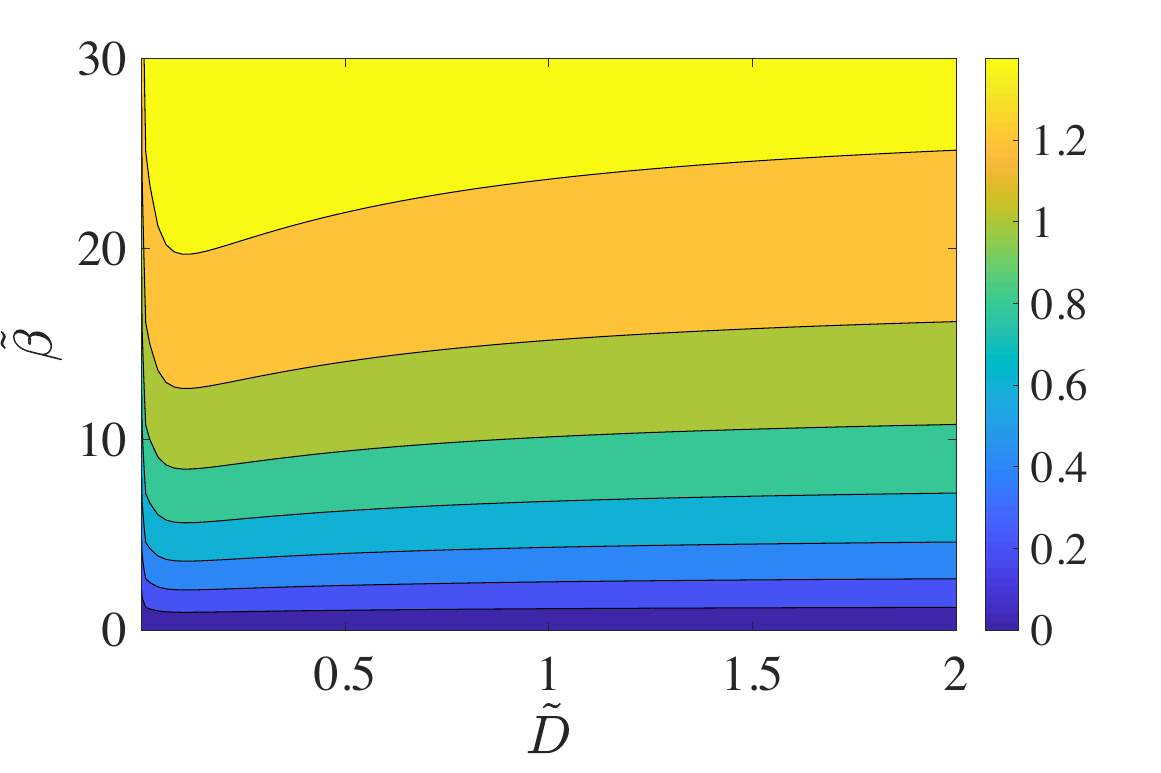}
\caption{case \ref{A3} -- level curves of $R_{0}$ as a function of $(\tilde D,\tilde\beta)$ for $\tilde\mu=50$ (values $R_{0,N}$ computed with $N=300$).}
\label{f__R0_vs_Dbeta}
\end{figure}

We can also combine the above results in a single plot, Figure \ref{f__R0_vs_Dbeta}, showing the level curves of the surface $R_{0}=R_{0}(\tilde D,\tilde\beta)$. The computation was performed on a grid of size $105\times100$, and took just few minutes (say less than $10$).

\begin{figure}[!ht]
\centering
\includegraphics[width=.48\textwidth]{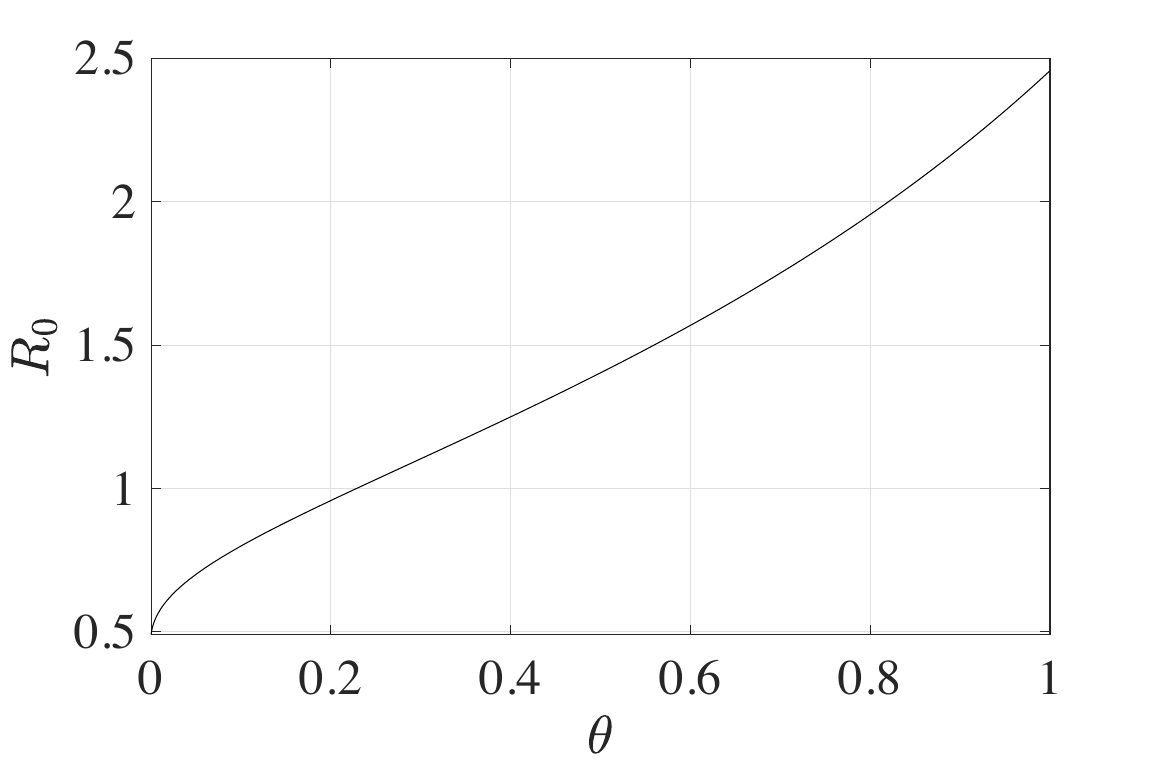}
\caption{case \ref{B3} -- $R_{0}$ as a function of $\theta$ (values $R_{0,N}$ computed with $N=200$).}
\label{f_theta_vs_R0}
\end{figure}

As for case \ref{B3}, in Figure \ref{f_theta_vs_R0} we show the variation of $R_{0}$ as a function of the vertical transmission probability $\theta$, confirming the biological expectation that this function is monotonically increasing as it is explicitly for case \ref{B1}. Notice that an increase of the vertical transmission can lead to the spread of the infection.
The grid size is $P=1000$ again, and the computation took less than $2$ minutes.

\begin{figure}[!ht]
\centering
\includegraphics[width=.48\textwidth]{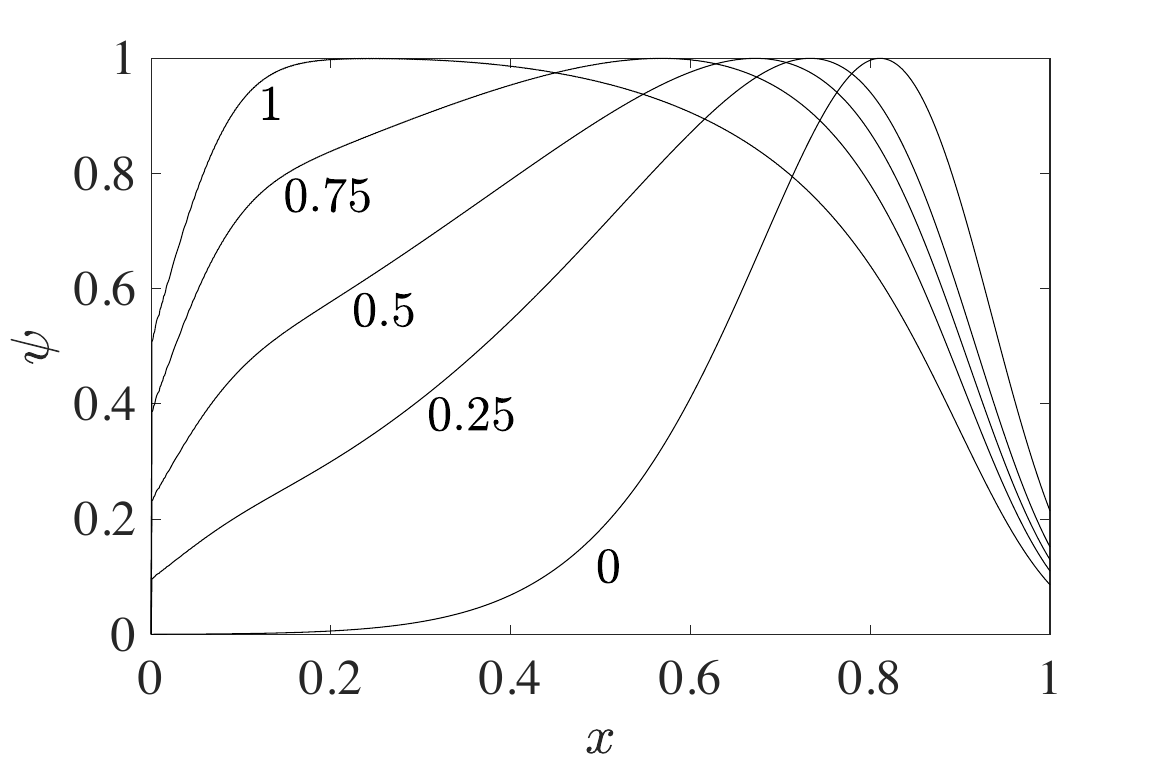}
\includegraphics[width=.48\textwidth]{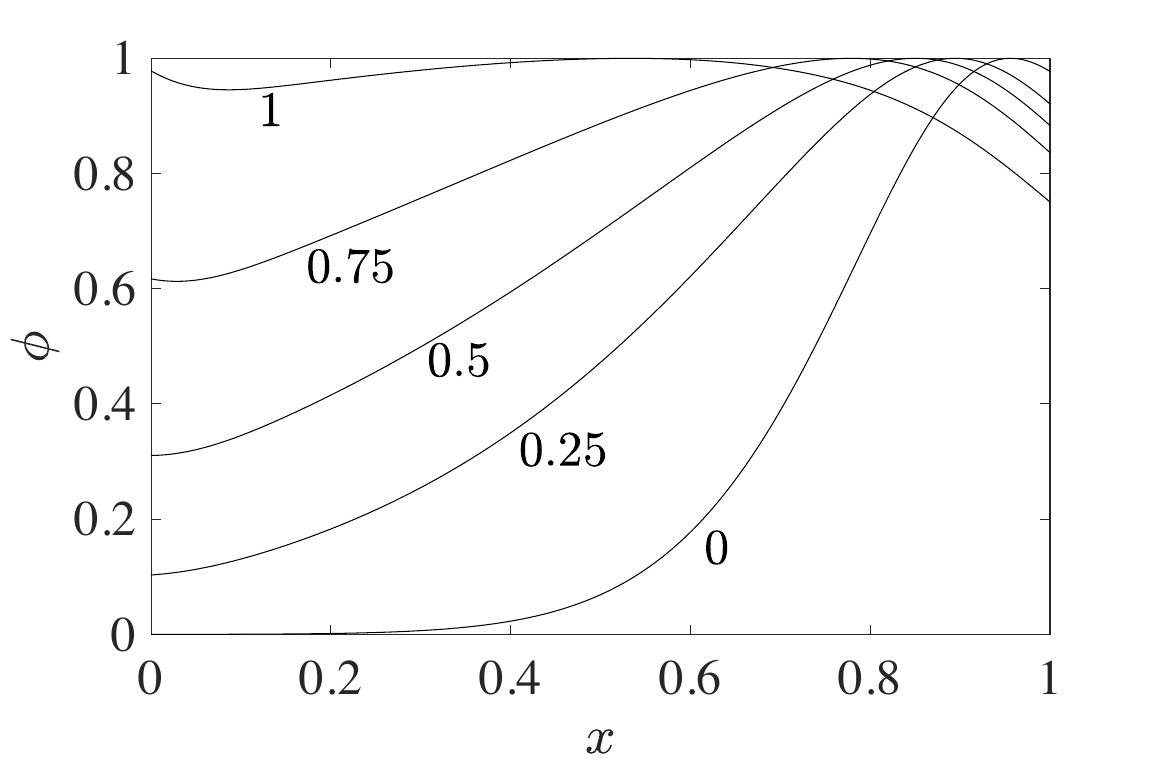}
\caption{case \ref{B3} -- normalized eigenfunctions (left) and generalized eigenfunctions (right) for several values of $\theta$, computed with $N=1\,000$.}
\label{f_B_eigenfunctions}
\end{figure}

Finally, always for case \ref{B3}, we show in Figure \ref{f_B_eigenfunctions} both the positive eigenfunctions (left) and the positive generalized eigenfunctions (right), normalized to have (absolute) maximum $1$. The vertical transmission probability is varied in the set $\{0,0.25,0.5,0.75,1\}$, but independently of $\theta$ we can appreciate how the proposed approach approximate smooth curves as expected.
%-----------------------------------------------------------------------------
%-----------------------------------------------------------------------------

\smallskip
\noindent{\bf Acknowledgement:} DB and RV are members of INdAM Research group GNCS and are partially supported by the INdAM GNCS project ``Approssimazione numerica di problemi di evoluzione: aspetti deterministici e stocastici'' (2018). JR is part of the Catalan research group 2017 SGR 1392 and has partially received support from the Spanish Ministry of Science and Innovation, reference MTM2017-84214-C2, and from the University of Girona, references MPC UdG 2016/047 and PONT2019/08.
%-----------------------------------------------------------------------------
%-----------------------------------------------------------------------------
%-----------------------------------------------------------------------------
\bibliographystyle{abbrv}
\bibliography{R0bib}
%-----------------------------------------------------------------------------
\end{document}